\documentclass[british,aps,manuscript,aps,rmp,preprint,tightenlines]{revtex4}
\usepackage[T1]{fontenc}
\usepackage[latin9]{inputenc}
\usepackage{units}
\usepackage{amsthm}
\usepackage{amsmath}
\usepackage{amssymb}
\usepackage{wasysym}
\usepackage[pdftex]{graphicx}
\usepackage{esint}

\makeatletter

\newcommand{\noun}[1]{\textsc{#1}}

\@ifundefined{textcolor}{}
{%
 \definecolor{BLACK}{gray}{0}
 \definecolor{WHITE}{gray}{1}
 \definecolor{RED}{rgb}{1,0,0}
 \definecolor{GREEN}{rgb}{0,1,0}
 \definecolor{BLUE}{rgb}{0,0,1}
 \definecolor{CYAN}{cmyk}{1,0,0,0}
 \definecolor{MAGENTA}{cmyk}{0,1,0,0}
 \definecolor{YELLOW}{cmyk}{0,0,1,0}
}
\numberwithin{equation}{section}
\numberwithin{figure}{section}
\usepackage{enumitem}		

\makeatother

\usepackage{babel}
\begin{document}

\title{Modelling elastic structures with strong nonlinearities with application
to stick-slip friction}

\author{Robert Szalai}

\affiliation{University of Bristol, Queen's Bldg., University Walk, Bristol, BS8
1TR, UK, email: }

\email{r.szalai@bristol.ac.uk}

\date{5th September 2013}
\begin{abstract}
An exact transformation method is introduced that reduces the governing
equations of a continuum structure coupled to strong nonlinearities
to a low dimensional equation with memory. The method is general and
well suited to problems with point discontinuities such as friction
and impact at point contact. It is assumed that the structure is composed
of two parts: a continuum but linear structure and finitely many discrete
but strong nonlinearites acting at various contact points of the elastic
structure. The localised nonlinearities include discontinuities, e.g.,
the Coulomb friction law. Despite the discontinuities in the model,
we demonstrate that contact forces are Lipschitz continuous in time
at the onset of sticking for certain classes of structures. The general
formalism is illustrated for a continuum elastic body coupled to a
Coulomb-like friction model.
\end{abstract}
\maketitle

\section{Background}

One of the greatest concerns of engineers is modelling friction and
impact. These two strong nonlinearities occur in many mechanical structures,
e.g., underplatform dampers of turbine blades \citep{Firrone,PetrovFEA},
tyre models \citep{PacejkaMagic}, or in general any jointed structure
\citep{Quinn,Segalman}. The most common way of modelling such systems
is to take a finite dimensional approximation assuming that the omitted
dynamics has only a small effect on the overall result. Such models
can then be analysed using the well established theory of non-smooth
dynamical systems \citep{diBernardoBook,FilippovBook}. The application
of this theory to engineering structures provides a great insight,
even though not all phenomena could be experimentally confirmed \citep{Popp}.
Recent results however indicate significant deficiencies that question
the predictive power of finite dimensional non-smooth models. It was
shown that when a rigid constraint becomes slightly compliant in a
friction-type system small-scale instabilities develop \citep{SieberSingular}.
This means that refining the model by including more degrees of freedom
can lead to qualitatively disagreeing results. The solution can also
become non-deterministic \citep{ColomboJeffrey,PainleveAlan} or non-unique
in forward time for an otherwise well specified initial condition.
Therefore a better modelling framework is necessary that either eliminates
inconsistency and non-determinism or at least provides a hint about
the physical mechanism that causes such behaviour.

The most apparent problem with finite dimensional non-smooth models
of mechanical systems is that they use rigid body dynamics to describe
the motion. This includes finite mode approximation of elastic structures,
where each mode has a non-zero modal mass \citep{EwinsBook}. When
two contacting elastic bodies slip and then suddenly stick their contact
points will experience a jump in acceleration. In case of a finite
mass at the contact point, the contact force also has a jump. In reality,
however, the mass of the contact point or contact surface is zero,
which implies that the contact force must be continuous. For this
reason standard finite mode description of elastic bodies is qualitatively
inaccurate. The continuity of contact force should be preserved by
mechanical models.

In this paper we propose a formalism that helps better understand
and perhaps solve the above problems. We investigate mechanical systems
that consist of linear elastic structures coupled at isolated points
of contact with strong nonlinearities. This class of problems include
mechanisms with Coulomb-like friction models. The dynamics of impact
is considered in a companion paper \citep{SzalaiMZImpact}. Our formalism
accounts for the zero mass of the contact point without artificially
introducing coupling springs as in \citep{Melcher} to regularise
the problem. To achieve such model reduction and still provide an
exact description we introduce memory terms. Within this new framework
the dynamics is described by a low-dimensional delay equation. We
show that in our formulation contact forces are Lipschitz or continuous
for certain classes of structures when the solution transitions onto
the discontinuity surface. The new formalism also leads to well defined
dynamics, since small perturbations of the reduced model do not affect
the qualitative features of the dynamics in general.

Time-delayed models have already been in use when modelling friction
\citep{HessSoom,Putelat2011}. In these cases, however the delay parameters
are fitted to experimental observations. We hope that through our
theory these empirical parameters can gain physical meaning.

The paper is organised as follows. In section \ref{sec:Mechmod} we
present our general mechanical model. In section \ref{sec:ModRed}
we describe our model reduction technique, discuss the convergence
of the method and its implications to non-smooth systems. The derivation
of the memory term is illustrated through the examples of a pre-tensed
string and a cantilever beam. In section \ref{sec:Bowstring} we present
the example of a bowed string. We demonstrate the properties of the
transformed equation of motion in particular its convergence as the
number of vibration modes goes to infinity.

\section{Mechanical model\label{sec:Mechmod}}

\begin{figure}
\begin{centering}
\includegraphics[scale=0.8]{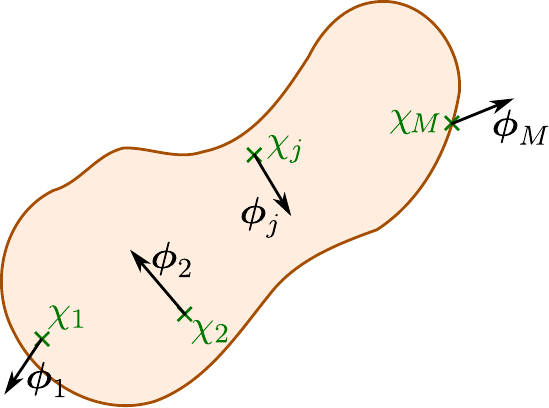}
\par\end{centering}

\caption{(colour online) A linear elastic structure with contact points $\chi_{j}$.
Each contact point has a three-dimensional motion, however we project
that motion to vectors $\boldsymbol{\phi}_{j}$ to obtain a scalar
valued resolved variable $y_{j}$. If we need to resolve more than
one direction of the motion of the contact point $\chi_{j}$, we attach
multiple labels $\chi_{j}=\chi_{j+1}=\chi_{j+2}$ to the same point.
Then the motion is projected by the linearly independent vectors $\phi_{j},\,\phi_{j+1},\,\phi_{j+2}$
to yield the resolved variables $y_{j},\, y_{j+1},\, y_{j+2}$. By
this definition we ensure a one-to-one mapping between indices of
labels $\chi_{j}$ and vectors $\phi_{j}$.\label{fig:potato}}
\end{figure}
The mechanical model of our structure is divided into two parts, a
linear elastic body and several discrete non-smooth nonlinearities
that are coupled to the continuum structure. First we describe our
assumptions on the linear but infinite dimensional part of the model
and then we explain how non-smooth nonlinearities are coupled to the
system. The description is sufficiently general to describe friction
oscillators and impact phenomena. For simplicity, we only focus on
a single elastic structure, but our framework is trivially extensible
to mechanisms involving multiple linear structures coupled at (strongly)
nonlinear joints.

We assume that the displacement of a material point $\chi$ of the
structure at time $t$ is represented by $\boldsymbol{u}(\chi,t)$.
We also assume that the motion $\boldsymbol{u}(\chi,t)$ can be expressed
as a series
\begin{equation}
\boldsymbol{u}(\chi,t)=\sum_{k=1}^{\infty}\boldsymbol{\psi}_{k}(\chi)x_{k}(t),\label{eq:dispField}
\end{equation}
where $\boldsymbol{\psi}_{k}(\chi)$ are three dimensional vector
valued functions depending on the spatial coordinates of the structure
only. The \emph{generalised coordinates} $x_{k}$ can be arranged
into a vector $\boldsymbol{x}=(x_{1},x_{2},\ldots)^{T}\in\mathbb{R}^{\infty}$
to simplify the notation. Due to linearity the governing equations
can be written as
\begin{equation}
\ddot{\boldsymbol{x}}(t)+\boldsymbol{C}\dot{\boldsymbol{x}}(t)+\boldsymbol{K}\boldsymbol{x}(t)=\boldsymbol{f}_{e}(t),\label{eq:GenLinSystem}
\end{equation}
where $\boldsymbol{C}$ and $\boldsymbol{K}$ are the damping and
stiffness matrices, respectively, both assumed being multiplied by
the inverse mass matrix from the left. The forcing term $\boldsymbol{f}_{e}(t)$
acts as a placeholder for the non-smooth part of the system and will
be replaced with with specific terms. Equation (\ref{eq:GenLinSystem})
allows for internal resonances. We assume that these resonances are
restricted to arbirarily large but finite dimensional subspaces of
the state space, which is necessary to guarantee basic convergence
properties of the solution $ $$\boldsymbol{x}(t)$ as shown in appendix
\ref{sub:C0semigroup}. We also assume that (\ref{eq:GenLinSystem})
is stable in the Lyapunov sense for the same reason.

When $\boldsymbol{C}$ and $\boldsymbol{K}$ matrices are simultaneously
diagonalisable the equation of motion can be written in the form of
\begin{equation}
\ddot{\boldsymbol{x}}(t)+2\boldsymbol{D}\boldsymbol{\Omega}\dot{\boldsymbol{x}}(t)+\boldsymbol{\Omega}^{2}\boldsymbol{x}(t)=\boldsymbol{f}_{e}(t),\label{eq:linModeDescr}
\end{equation}
where $\boldsymbol{\Omega}=\mathrm{diag}(\omega_{1},\omega_{2},\ldots)$
and $\boldsymbol{D}=\mathrm{diag}(D_{1},D_{2},\ldots)$. In the unforced
case ($\boldsymbol{f}_{e}=\boldsymbol{0}$), the vector components
of $\boldsymbol{x}$ on the left side of equation (\ref{eq:linModeDescr})
are decoupled, which means that the homogeneous equation can be solved
for each $x_{k}$ independently. Therefore $x_{k}(t)$ and $\boldsymbol{\psi}_{k}(\chi)$
are called the \emph{mode}s and \emph{mode shape}s of the system,
respectively, with $\omega_{k}$ natural frequencies and $D_{k}$
damping ratios \citep{EwinsBook}. System (\ref{eq:linModeDescr})
are called \emph{modal equations} describing the motion through the
\emph{modal coordinates} $\boldsymbol{x}$. Our results do not require
that the equations of motion assume the form (\ref{eq:linModeDescr}),
however in section \ref{sub:L1conv} some restriction on the eigenvalues
of (\ref{eq:GenLinSystem}) is necessary to characterise the convergence
of the reduced equations of motion.

In order to take into account the coupling of the contact points to
non-smooth nonlinearities we need to characterise their motion. For
simplicity we assume point contact. Let us denote the motion of the
$j$-th contact point at $\chi_{j}$ along the direction of vector
$\boldsymbol{\phi}_{j}$ by 
\begin{equation}
y_{j}(t)=\boldsymbol{\phi}_{j}\cdot\boldsymbol{u}(\chi_{j},t)\label{eq:resolvDef}
\end{equation}
as illustrated in Fig. \ref{fig:potato}. We call the positions $y_{j}(t)$
and the velocities $\dot{y}_{j}(t)$ of the contact points \emph{resolved
variables}. We assume $M$ contact points, thus we define $\boldsymbol{y}=(y_{1},\ldots,y_{M},\dot{y}_{1},\ldots,\dot{y}_{M})^{T}$.
Substituting (\ref{eq:dispField}) into (\ref{eq:resolvDef}) we obtain
the motion of the contact points as a function of the solution of
equation (\ref{eq:GenLinSystem}) 
\begin{equation}
y_{j}(t)=\boldsymbol{n}{}_{j}\cdot\boldsymbol{x}(t),\label{eq:resolvDotProd}
\end{equation}
where
\begin{equation}
\boldsymbol{n}{}_{j}=\left(\boldsymbol{\phi}_{j}\cdot\boldsymbol{\psi}_{1}(\chi_{j}),\boldsymbol{\phi}_{j}\cdot\boldsymbol{\psi}_{2}(\chi_{j}),\ldots\right)^{T}.
\end{equation}
Vectors $\boldsymbol{n}{}_{j}$ are assumed to be linearly independent,
spanning an $M$ dimensional subspace of $\mathbb{R}^{\infty}$.

The nonlinearities are incorporated into the model through the forcing
term $\boldsymbol{f}_{e}(t)$. We assume that the nonlinearities only
depend on the resolved variables. They are also piecewise continuous
with isolated discontinuities. Thus the contact forces acting at contact
points $\chi_{j}$ in the direction $\boldsymbol{\phi}_{j}$ are written
as $f_{j}(\boldsymbol{y}(t),t).$ Summing up all the nonlinearities
completes our model description by providing the right-hand side of
equation (\ref{eq:linModeDescr}) in the form of 
\begin{equation}
\boldsymbol{f}_{e}(t)=\sum_{j=1}^{M}\boldsymbol{n}_{j}f_{j}(\boldsymbol{y}(t),t).\label{eq:interForce}
\end{equation}

\section{Reduction of the mechanical model\label{sec:ModRed}}

We aim to reduce the number of dimensions of our mechanical model
(\ref{eq:interForce}) to an equation that only involves the $2M$
number of resolved variables all contained in the vector $\boldsymbol{y}$.
To achieve this we use the Mori-Zwanzig formalism \citep{ChorinPNAS}
to arrive at a $2M$ dimensional first order delay equation. The solution
of the reduced system agrees with the solution of the full system
for the resolved variables, while the rest of the variables are discarded.
Our method can be viewed as a way of producing a Green's function
for only parts of the system. In this formalism a convolution with
the resolved variables represents the effect of the eliminated variables
on the dynamics of the resolved variables. The technical details of
the transformation are described in appendices \ref{sec:AppTrafo}
and \ref{sec:AppMemKer}.

To simplify our calculation we transform equation (\ref{eq:linModeDescr})
into a first order form of 
\begin{equation}
\dot{\boldsymbol{z}}(t)=\boldsymbol{R}\boldsymbol{z}(t)+\boldsymbol{f}(t),\label{eq:1storderODE}
\end{equation}
where
\begin{equation}
\boldsymbol{R}=\left(\begin{array}{cc}
\boldsymbol{0} & \boldsymbol{I}\\
-\boldsymbol{K} & -\boldsymbol{C}
\end{array}\right)\;\mbox{and}\;\boldsymbol{f}(t)=\left(\begin{array}{c}
\boldsymbol{0}\\
\boldsymbol{f}_{e}(t)
\end{array}\right).\label{eq:RmatDef}
\end{equation}

We already have a way of obtaining the resolved variables from the
generalised coordinates $\boldsymbol{x}$ through a dot product with
vectors $\boldsymbol{n}{}_{j}$ as shown in equation (\ref{eq:resolvDotProd}).
In a matrix-vector notation we write the conversion as 
\begin{equation}
\boldsymbol{y}(t)=\boldsymbol{V}\boldsymbol{z}(t),\;\mbox{with}\;\boldsymbol{V}=\left(\begin{array}{c}
\boldsymbol{v}_{1}^{T}\\
\vdots\\
\boldsymbol{v}_{2M}^{T}
\end{array}\right)\,\mbox{and}\,\boldsymbol{v}_{j}=\left(\begin{array}{c}
\boldsymbol{n}_{j}\\
\boldsymbol{0}
\end{array}\right),\;\boldsymbol{v}_{M+j}=\left(\begin{array}{c}
\boldsymbol{0}\\
\boldsymbol{n}_{j}
\end{array}\right).
\end{equation}
To obtain our reduced model we construct a projection matrix $\boldsymbol{S}$
that acts on the generalised coordinates and has a $2M$ dimensional
range. In order to do that we also define a lifting operator in the
form of
\begin{equation}
\boldsymbol{z}(t)=\boldsymbol{W}\boldsymbol{y}(t),\;\mbox{with\;}\boldsymbol{W}=\left(\begin{array}{ccc}
\boldsymbol{w}_{1} & \cdots & \boldsymbol{w}_{2M}\end{array}\right).
\end{equation}
Note that the lifting operator with $\boldsymbol{W}$ does not reproduce
the full solution from the resolved variables, it is merely used as
a technical tool. Moreover, $\boldsymbol{W}$ does not depend on the
physical system, it can be chosen to suite the reduction procedure.
By combining matrices $\boldsymbol{V}$ and $\boldsymbol{W}$ we obtain
our projections $\boldsymbol{S}=\boldsymbol{W}\boldsymbol{V}$ and
$\boldsymbol{Q}=\boldsymbol{I}-\boldsymbol{S}$ on the condition that
$\boldsymbol{m}_{k}$ satisfy $\boldsymbol{m}_{k}\cdot\boldsymbol{n}_{l}=0$
if $k\neq l$ and $\boldsymbol{m}_{k}\cdot\boldsymbol{n}_{k}=1$.
This constraint can also be expressed as $\boldsymbol{I}=\boldsymbol{V}\boldsymbol{W}$
(note the order of the two matrices).

To guarantee that the terms in the reduced equation are bounded the
choice of $\boldsymbol{W}$ needs to be further restricted. Therefore
we assume that the range of $\boldsymbol{W}$ is invariant under $\boldsymbol{R}$,
that is, 

\begin{equation}
\boldsymbol{R}\boldsymbol{w}_{j}\in\mathrm{span}(\boldsymbol{w}_{1},\boldsymbol{w}_{2},\ldots,\boldsymbol{w}_{2M}),\;\mbox{for}\; j=1,2,\ldots,2M.\label{eq:RangeCond}
\end{equation}
Equivalently, $\boldsymbol{w}_{j}$ can be constructed as linear combinations
of $2M$ number of eigenvectors of $\boldsymbol{R}$, since eigenvectors
are invariant by definition. This assumption is key to our analysis. 

In case of the modal equations (\ref{eq:linModeDescr}) the columns
of $\boldsymbol{W}$ can be explicitly constructed in the form of
\begin{equation}
\boldsymbol{w}_{j}=\left(\begin{array}{c}
\boldsymbol{m}_{j}\\
\boldsymbol{0}
\end{array}\right)\;\mbox{and}\;\boldsymbol{w}_{M+j}=\left(\begin{array}{c}
\boldsymbol{0}\\
\boldsymbol{m}_{j}
\end{array}\right).\label{eq:WmDef}
\end{equation}
Due to the block diagonal form of $\boldsymbol{R}$ condition (\ref{eq:RangeCond})
holds when the $\boldsymbol{m}_{j}$ vectors are chosen such that
they only have $M$ number of non-zero components: 
\begin{equation}
\left[\boldsymbol{m}_{j}\right]_{p}=0,\quad\mbox{for}\quad p<P\;\mbox{and}\; p\ge P+M.\label{eq:mConstraint}
\end{equation}

As the last step before arriving at the reduced equations we define
\begin{align}
\boldsymbol{A} & =\boldsymbol{V}\boldsymbol{R}\boldsymbol{W}, &  & \in\mathbb{R}^{2M\times2M}\label{eq:Amat-1}\\
\boldsymbol{H}(t)\boldsymbol{z}(s) & =\boldsymbol{V}\boldsymbol{R}\boldsymbol{Q}\mathrm{e}^{\boldsymbol{R}(t-s)}\boldsymbol{z}(s), &  & \in\mathbb{R}^{2M}\label{eq:noise}\\
\boldsymbol{L}_{j}(\tau) & =\boldsymbol{V}\mathrm{e}^{\boldsymbol{R}\tau}\boldsymbol{v}_{M+j}-\boldsymbol{A}\boldsymbol{V}\int_{0}^{\tau}\mathrm{e}^{\boldsymbol{R}\theta}\boldsymbol{v}_{M+j}\mathrm{d}\theta, &  & \in\mathbb{R}^{2M}\label{eq:forceKern}
\end{align}
where $\boldsymbol{z}(s)$ is the initial condition of equation (\ref{eq:1storderODE})
at time $s$ and $\mathrm{e}^{\boldsymbol{R}t}$ is the fundamental
matrix of equation (\ref{eq:1storderODE}). The matrix $\boldsymbol{A}$
is bounded if condition (\ref{eq:RangeCond}) holds, while $\boldsymbol{H}(t)\boldsymbol{z}(s)$
is bounded if the initial condition $\boldsymbol{z}(s)$ is bounded,
too. To obtain the memory kernel $\boldsymbol{L}_{j}(\tau)$, one
needs to solve the first order system (\ref{eq:1storderODE}) for
$M$ different initial conditions $\boldsymbol{v}_{M+j}$. This solution
may not be bounded which we rectify by integrating it (see sections
\ref{sub:LjConst} and \ref{sub:L1conv}).

Using the expression (\ref{eq:interForce}) of the forcing term our
reduced equation that is equivalent to (\ref{eq:1storderODE}) in
the resolved variables becomes
\begin{equation}
\dot{\boldsymbol{y}}(t)=\boldsymbol{A}\boldsymbol{y}(t)+\sum_{j=1}^{M}\left[\boldsymbol{L}_{j}(0)f_{j}(\boldsymbol{y}(t),t)+\int_{0}^{t-s}\mathrm{d}_{\tau}\boldsymbol{L}_{j}(\tau)f_{j}(\boldsymbol{y}(t-\tau),t-\tau)\right]+\boldsymbol{H}(t)\boldsymbol{z}(s),\label{eq:MZequation}
\end{equation}
where the integral is understood in the Riemann-Stieltjes sense (see
section \ref{sub:Stieltjes}). The formal equivalence of (\ref{eq:1storderODE})
and (\ref{eq:MZequation}) is proved in appendix \ref{sec:AppTrafo}
and the formulae of $\boldsymbol{A}$, $\boldsymbol{H}(t)\boldsymbol{z}(s)$
and $\boldsymbol{L}_{j}(\tau)$ are derived in appendix \ref{sec:AppMemKer}.
It turns out that for mechanical systems the alternative form of (\ref{eq:MZequation})
given below by equation (\ref{eq:MZL1eq}) is more appropriate.

\subsection{The meaning of the Riemann-Stieltjes integral\label{sub:Stieltjes}}

To provide some intuition about the interpretation of the integral
in (\ref{eq:MZequation}) we remark that if $\boldsymbol{L}_{j}(\tau)$
is differentiable, the symbol $\mathrm{d}_{\tau}\boldsymbol{L}_{j}(\tau)$
can be replaced by $\nicefrac{\mathrm{d}}{\mathrm{d}\tau}\boldsymbol{L}_{j}(\tau)\cdots\mathrm{d}\tau$
to arrive at an ordinary Riemann integral. Similar simplification
can be achieved if $\boldsymbol{L}_{j}(\tau)$ is Lipschitz continuous,
however $\mathrm{d}_{\tau}\boldsymbol{L}_{j}(\tau)$ is replaced by
the left derivative of $\boldsymbol{L}_{j}(\tau)$. A discontinuity
of $\boldsymbol{L}_{j}(\tau)$ at $\tau^{\star}$ translates into
a discrete delay term $\mathfrak{L}_{j}(\tau^{\star})f_{j}(\boldsymbol{y}(t-\tau^{\star}),t-\tau^{\star})$,
where $\mathfrak{L}_{j}(\tau^{\star})$ is the gap in the function
at $\tau^{\star}$. Therefore if $\boldsymbol{L}_{j}(\tau)$ is piecewise
differentiable (or piecewise Lipschitz), with finite number of isolated
discontinuities, the integral in (\ref{eq:MZequation}) can be replaced
by a sum of Riemann integrals and discrete delays. In particular,
assuming that the discontinuities occur at $\tau_{p}^{\star}$, $p=1,\ldots,P$
and that $\tau_{0}^{\star}=0$ and $\tau_{P+1}^{\star}=t$ we can
write 
\begin{equation}
\int_{0}^{t-s}\mathrm{d}_{\tau}\boldsymbol{L}_{j}(\tau)f(t-\tau)\mathrm{d}\tau=\sum_{p=0}^{P}\int_{\tau_{p}^{\star}}^{\tau_{p+1}^{*}-s}\left[\frac{\mathrm{d}}{\mathrm{d}\tau}\boldsymbol{L}_{j}(\tau)\right]f(t-\tau)\mathrm{d}\tau+\sum_{p=1}^{P}\mathfrak{L}_{j}(\tau_{p}^{\star})f(t-\tau_{p}^{\star}),
\end{equation}
where $\mathfrak{L}_{j}(\tau^{\star})=\lim_{\tau\to\tau*+0}\boldsymbol{L}_{j}(\tau)-\lim_{\tau\to\tau*-0}\boldsymbol{L}_{j}(\tau)$.

One can interpret the memory kernels $\boldsymbol{L}_{j}(\tau)$ as
an indication how waves travel within the structure from contact point
$\chi_{j}$ to all contact points including them returning to $\chi_{j}$.
Indeed, the tail part of $\boldsymbol{L}_{j}(\tau)$ for $\tau>0$
describes how a force applied in the past is affecting the contact
points at current time. This notion is therefore analogous to having
waves that depart from point $\chi_{j}$ and arrive at all contact
points $\chi_{l}$, $l=1,\ldots,M$ $\tau$ time later. In particular
if waves do not disperse, they arrive simultaneously at a give material
point $\chi_{l}$ exactly $\tau^{\star}$ time later and that corresponds
to a discrete delay represented by a discontinuity in $\boldsymbol{L}_{j}(\tau)$
at $\tau=\tau^{\star}$. 

The time delay in (\ref{eq:MZequation}) tends to infinity, since
the history of $\boldsymbol{y}(t)$ to be taken into account grows
with time. However, if $\boldsymbol{L}_{j}(\tau)$ tends to a constant
as $\tau\to\infty$, the integrals in (\ref{eq:MZequation}) may be
truncated at a finite delay time to produce an approximate model.
In case of a conservative equation (\ref{eq:linModeDescr}), the delay
cannot be truncated at a finite time, because the motion of the free
structure will never stop due to the conservation of the kinetic energy.
This is illustrated in section \ref{sub:StringVib}.

\subsection{An alternative form of the reduced equations\label{sub:LjConst}}

It was mentioned before that $\boldsymbol{L}_{j}(\tau)$ may not be
a bounded function. Therefore we integrate the convolution in equation
(\ref{eq:MZequation}) by parts, so that the integral of $\boldsymbol{L}_{j}(\tau)$
and the derivative of the contact force $f_{j}$ appear in the rewritten
equation. However, to be able to integrate $\boldsymbol{L}_{j}(\tau)$
and ensure that the integral does not grow out of bound, we need to
decompose $\boldsymbol{L}_{j}(\tau)$ into a constant and an oscillatory
term 
\begin{equation}
\boldsymbol{L}_{j}(\tau)=\boldsymbol{L}_{j}^{\infty}+\boldsymbol{L}_{j}^{0}(\tau),\label{eq:L0Def}
\end{equation}
where
\begin{equation}
\boldsymbol{L}_{j}^{\infty}=\lim_{t\to\infty}\frac{1}{t}\int_{0}^{t}\boldsymbol{L}_{j}(\tau)\mathrm{d}\tau.\label{eq:LinfLimit}
\end{equation}
Since we have a formal expression for $\boldsymbol{L}_{j}(\tau)$
in the form of equation (\ref{eq:forceKern}), we can calculate the
integral
\begin{equation}
\int_{0}^{t}\boldsymbol{L}_{j}(\tau)\mathrm{d}\tau=\int_{0}^{t}\left(\boldsymbol{V}\mathrm{e}^{\boldsymbol{R}\tau}-\boldsymbol{A}\boldsymbol{V}\int_{0}^{\tau}\mathrm{e}^{\boldsymbol{R}\theta}\mathrm{d}\theta\right)\mathrm{d}\tau\boldsymbol{v}_{M+j}.\label{eq:LjIntegral}
\end{equation}
From the simple rule of differentiating an exponential we find that
$\int_{0}^{\tau}\mathrm{e}^{\boldsymbol{R}\theta}\mathrm{d}\theta=\boldsymbol{R}^{-1}\left(\mathrm{e}^{\boldsymbol{R}\tau}-\boldsymbol{I}\right)$.
Using this formula twice in equation (\ref{eq:LjIntegral}), we are
left with
\begin{equation}
\int_{0}^{t}\boldsymbol{L}_{j}(\tau)\mathrm{d}\tau=\left(\boldsymbol{V}\boldsymbol{R}^{-1}\left(\mathrm{e}^{\boldsymbol{R}t}-\boldsymbol{I}\right)-\boldsymbol{A}\boldsymbol{V}\boldsymbol{R}^{-1}\left(\boldsymbol{R}^{-1}\left(\mathrm{e}^{\boldsymbol{R}t}-\boldsymbol{I}\right)-t\right)\right)\boldsymbol{v}_{M+j}.\label{eq:LjIntegralFull}
\end{equation}
Assuming that $\boldsymbol{V}\int_{0}^{\tau}\mathrm{e}^{\boldsymbol{R}\theta}\mathrm{d}\theta\boldsymbol{v}_{M+j}$
is bounded (see section \ref{sub:L1conv}) the limit (\ref{eq:LinfLimit})
becomes
\begin{equation}
\boldsymbol{L}_{j}^{\infty}=\boldsymbol{A}\boldsymbol{V}\boldsymbol{R}^{-1}\boldsymbol{v}_{M+j}.\label{eq:LinfDef}
\end{equation}
Note that $\boldsymbol{V}\boldsymbol{R}^{-1}\boldsymbol{v}_{M+j}$
is the static displacement of the contact points under a static unit
load at contact point $\chi_{j}$. Indeed, by expanding equation $\boldsymbol{R}\left(\boldsymbol{x},\boldsymbol{y}\right)^{T}=\left(\boldsymbol{0},\boldsymbol{n}_{j}\right)^{T}$
we get $\boldsymbol{K}\boldsymbol{x}=\boldsymbol{n}_{j}$ and $\boldsymbol{y}=\boldsymbol{0}$,
where $\boldsymbol{x}$ is the static displacement of the generalised
coordinates, hence, $\boldsymbol{V}\boldsymbol{x}=\boldsymbol{V}\boldsymbol{R}^{-1}\boldsymbol{v}_{M+j}$
is the static displacement of the contact points.

Our definition (\ref{eq:L0Def}) of $\boldsymbol{L}_{j}^{0}(\tau)$
implies that either $\boldsymbol{L}_{j}^{0}(\tau)\to\boldsymbol{0}$
exponentially if the system is dissipative or $\boldsymbol{L}_{j}^{0}(\tau)$
oscillates with zero mean. Therefore we define the integral 
\begin{equation}
\boldsymbol{L}_{j}^{1}(t)=\int_{0}^{t}\boldsymbol{L}_{j}^{0}(\tau)\mathrm{d}\tau.\label{eq:L1def}
\end{equation}
The boundedness and smoothness of (\ref{eq:L1def}) is discussed in
section \ref{sub:L1conv}. By virtue of (\ref{eq:L1def}) the forcing
terms in (\ref{eq:MZequation}) can be integrated by parts as follows
\begin{multline}
\boldsymbol{L}_{j}(0)\boldsymbol{f}(t)+\int_{0}^{t-s}\mathrm{d}_{\tau}\boldsymbol{L}_{j}(\tau)f_{j}(\boldsymbol{y}(t-\tau),t-\tau)=\\
=\boldsymbol{L}_{j}^{\infty}f_{j}(t)+\boldsymbol{L}_{j}^{0}(t-s)f_{j}(s)+\int_{0}^{t-s}\mathrm{d}_{\tau}\boldsymbol{L}_{j}^{1}(\tau)\frac{\mathrm{d}}{\mathrm{d}t}[f_{j}(\boldsymbol{y}(t-\tau),t-\tau)],
\end{multline}
which has only bounded terms if the derivative $\frac{\mathrm{d}}{\mathrm{d}t}[f_{j}(t)]$
exists and is bounded almost everywhere and condition (\ref{eq:ResolventCondition})
in section \ref{sub:L1conv} holds. Using this transformation the
governing equation becomes
\begin{equation}
\dot{\boldsymbol{y}}(t)=\boldsymbol{A}\boldsymbol{y}(t)+\sum_{j=1}^{M}\left[\boldsymbol{L}_{j}^{\infty}f_{j}(t)+\int_{0}^{t-s}\mathrm{d}_{\tau}\boldsymbol{L}_{j}^{1}(\tau)\frac{\mathrm{d}}{\mathrm{d}t}[f_{j}(\boldsymbol{y}(t-\tau),t-\tau)]\right]+\boldsymbol{g}(t),\label{eq:MZL1eq}
\end{equation}
where $\boldsymbol{g}(t)=\boldsymbol{H}(t)\boldsymbol{z}(s)+\sum_{j=1}^{M}\boldsymbol{L}_{j}^{0}(t-s)f_{j}(s).$

It is important to note that according to the theory of delay equations
\citep{hale1993introduction} if the terms in equation (\ref{eq:MZL1eq})
are slightly perturbed, the solution of (\ref{eq:MZL1eq}) changes
only slightly under general conditions. This is a clear advantage
over the finite dimensional description where perturbation generally
leads to qualitatively disagreeing solutions \citep{SieberSingular}.

\subsection{Vibrations of a pre-tensed string\label{sub:StringVib}}

\begin{figure}
\begin{centering}
\includegraphics[width=0.99\linewidth]{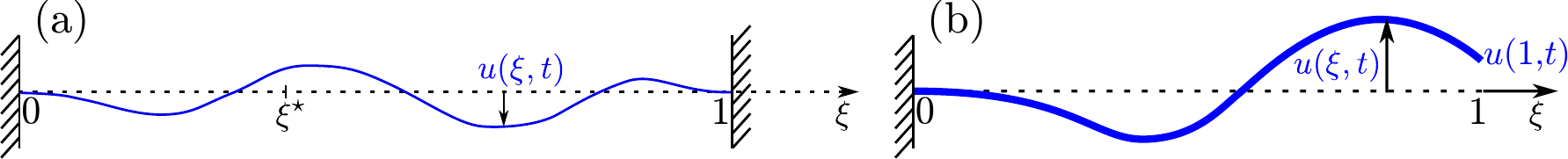}
\par\end{centering}

\caption{Schematic of a string (a) and a beam (b). The displacement of material
points is described by $u(\xi,t)$, which represents motion in the
direction of the arrows. \label{fig:StVenant}}

\end{figure}

In this section we illustrate the calculation of $\boldsymbol{A}$
and $\boldsymbol{L}_{j}^{1}(t)$ for a pre-tensed string without bending
stiffness. We keep the calculation as general as possible so that
it applies to systems with a single contact point that can be written
in the form of \ref{eq:linModeDescr}. The schematic of the string
is shown in Fig. \ref{fig:StVenant}(a), whose motion is described
by the equation 
\begin{equation}
\frac{\partial^{2}u}{\partial t^{2}}=c^{2}\frac{\partial^{2}u}{\partial\xi^{2}}+\delta(\xi-\xi^{\star})F_{c}(t)+\mbox{damping},\label{eq:StringPDE}
\end{equation}
where $c$ is the wave speed, $u(t,\xi)$ is the deflection of the
string, $t$ stands for time and $0\le\xi\le1$ is the coordinate
along the string. The boundary conditions $u(0,t)=0$ and $u(1,t)=0$
express that there is no movement at the two ends of the string. Equation
(\ref{eq:StringPDE}) indicates damping, which explicitly appears
in the mode decomposed system (\ref{eq:linModeDescr}) with non-zero
damping ratios. The string is forced at $\xi=\xi^{\star}$ by a contact
force $F_{c}(t)$, which is represented by the Dirac delta function
in equation (\ref{eq:StringPDE}).

The first step is to provide a mode decomposition in the form of equation
(\ref{eq:linModeDescr}), so that the vibration of the string is expressed
by (\ref{eq:dispField}), where the mode shapes are the scalar valued
$\psi_{k}(\xi)=\sin\left(k\pi\xi\right)$. The natural frequencies
of the system are $\omega_{k}=ck\pi$ and we assume uniform damping
$D_{k}=\nicefrac{1}{10}$ for all modes, which gives us
\begin{equation}
\boldsymbol{\Omega}=\mathrm{diag}(c\pi,c2\pi,\ldots),\quad\mbox{and}\quad\boldsymbol{D}=\mathrm{diag}(\nicefrac{1}{10},\nicefrac{1}{10},\ldots).
\end{equation}
The vibration at the contact point can be expressed as a linear combination
of all the modes $u(\xi^{\star},t)=\boldsymbol{n}\cdot\boldsymbol{x}(t)$,
where $\boldsymbol{n}=(\sin\left(\pi\xi^{\star}\right),\sin\left(2\pi\xi^{\star}\right),\ldots)^{T}$.
The resolved variables therefore are $y_{1}(t)=\boldsymbol{n}\cdot\boldsymbol{x}(t)$
and $y_{2}(t)=\dot{y}_{1}(t)$. We also choose $\boldsymbol{m}=(\nicefrac{1}{\sin\pi\xi^{\star}},0,0,\ldots)^{T}$
in formula (\ref{eq:WmDef}), thus 
\begin{equation}
\boldsymbol{A}=\left(\begin{array}{cc}
0 & 1\\
-\omega_{1}^{2} & -2D_{1}\omega_{1}
\end{array}\right).\label{eq:StringAmat}
\end{equation}

\begin{figure}
\begin{centering}
\includegraphics[width=0.99\linewidth]{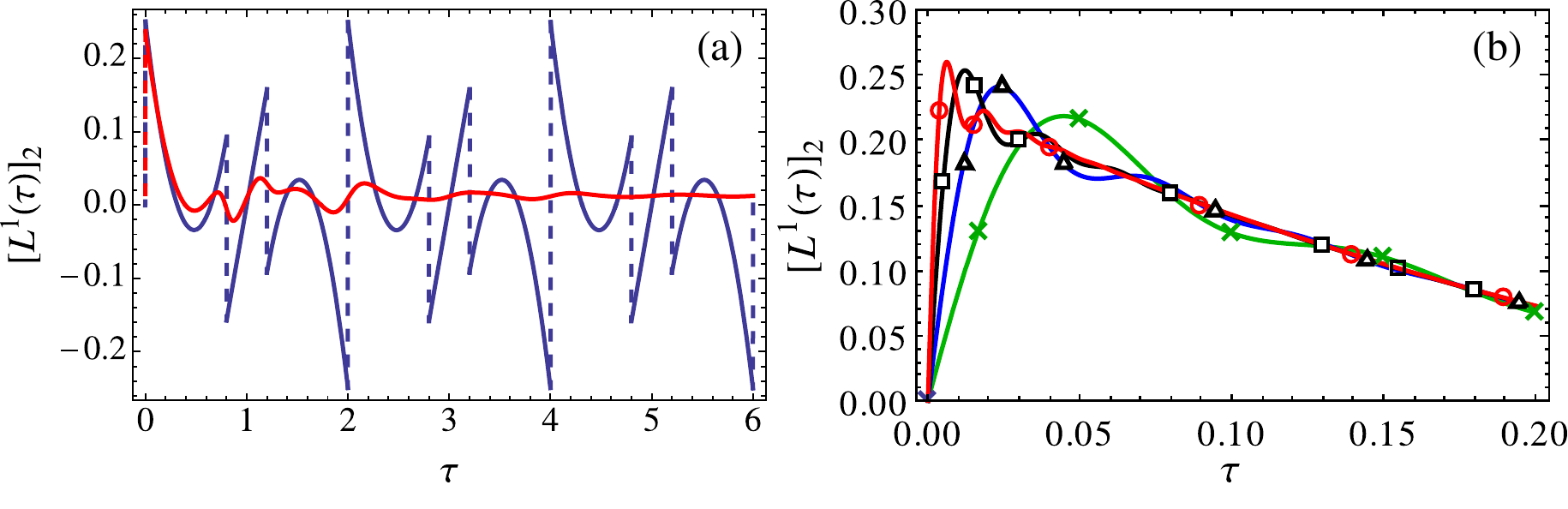}
\par\end{centering}

\caption{\label{fig:StringL1fig}(colour online) Graph of $\left[\boldsymbol{L}^{1}(\tau)\right]_{2}$
for the string equation (\ref{eq:StringPDE}) with at $\xi^{\star}=0.4,$
$c=1$. The blue line represents the conservative case and the red
corresponds to the damped case with $D_{k}=\nicefrac{1}{10}$. (b)
Graph of the function $\left[\boldsymbol{L}^{1}(\tau)\right]_{2}$
when truncating the series (\ref{eq:StringL2intGen}) at $20,40,80,160$
terms, denoted by $\times,\triangle,\square$ and $\ocircle$, respectively.}

\end{figure}

Next we calculate the function $\mathrm{e}^{\boldsymbol{R}t}\boldsymbol{v}_{2}$
that appears in the expression (\ref{eq:forceKern}) of $\boldsymbol{L}(\tau)$.
The equation whose solution we are seeking is $\dot{\boldsymbol{z}}=\boldsymbol{R}\boldsymbol{z}$,
which is expanded into
\begin{equation}
\dot{z}_{1,k}=z_{2,k},\quad\dot{z}_{2,k}=-2D_{k}\omega_{k}z_{2,k}-\omega_{k}^{2}z_{1,k},\quad k=1,2,3,\ldots.\label{eq:StringRExpand}
\end{equation}
The initial conditions that correspond to $\boldsymbol{z}(0)=\boldsymbol{v}_{2}$
are $z_{1,k}(0)=0$ and $z_{2,k}(0)=\left[\boldsymbol{n}\right]_{k}$.
The solution for the modes are
\begin{equation}
z_{1,k}(t)=\left[\boldsymbol{n}\right]_{k}\mathrm{e}^{-D_{k}\omega_{k}t}\frac{\sin\sqrt{1-D_{k}^{2}}\omega_{k}t}{\sqrt{1-D_{k}^{2}}\omega_{k}}.\label{eq:StringFullSol}
\end{equation}
Without assuming the form of $z_{1,k}(t)$ and evaluating formula
(\ref{eq:forceKern}) we find that
\begin{align}
\left[\boldsymbol{L}(t)\right]_{1} & =\boldsymbol{v}_{1}\cdot\mathrm{e}^{\boldsymbol{R}t}\boldsymbol{v}_{2}-\boldsymbol{v}_{2}\cdot\int_{0}^{t}\mathrm{e}^{\boldsymbol{R}\theta}\mathrm{d}\theta\boldsymbol{v}_{2}=0\quad\mbox{and}\nonumber \\
\left[\boldsymbol{L}(\tau)\right]_{2} & =\boldsymbol{v}_{2}\cdot\mathrm{e}^{\boldsymbol{R}t}\boldsymbol{v}_{2}+2D_{1}\omega_{1}\boldsymbol{v}_{1}\cdot\mathrm{e}^{\boldsymbol{R}t}\boldsymbol{v}_{2}+\omega_{1}^{2}\boldsymbol{v}_{1}\cdot\int_{0}^{t}\mathrm{e}^{\boldsymbol{R}\theta}\mathrm{d}\theta\boldsymbol{v}_{2}\nonumber \\
 & =\sum_{k=1}^{\infty}\left[\boldsymbol{n}\right]_{k}\left(\frac{d}{dt}z_{1,k}(t)+2D_{1}\omega_{1}z_{1,k}(t)+\omega_{1}^{2}\int_{0}^{t}z_{1,k}(t)d\theta\right).\label{eq:StringL2gen}
\end{align}
Eventually substituting (\ref{eq:StringFullSol}) into (\ref{eq:StringL2gen})
in the conservative case ($D_{k}=0$) $\left[\boldsymbol{L}(t)\right]_{2}$
of our system (\ref{eq:StringPDE}) becomes 
\begin{equation}
\left[\boldsymbol{L}(t)\right]_{2}=\sum_{k=1}^{\infty}\sin^{2}k\pi\xi^{\star}\frac{\omega_{1}^{2}}{\omega_{k}^{2}}+\sum_{k=2}^{\infty}\sin^{2}k\pi\xi^{\star}\left(1-\frac{\omega_{1}^{2}}{\omega_{k}^{2}}\right)\cos\omega_{k}t,\label{eq:StringL2}
\end{equation}
which is a divergent Fourier series, therefore equation (\ref{eq:MZequation})
cannot be utilised to describe the dynamics. The constant term in
equation (\ref{eq:StringL2gen}) regardless of the damping ratios
assumes the form 
\begin{equation}
\left[\boldsymbol{L}^{\infty}\right]_{2}=-\omega_{1}^{2}\boldsymbol{v}_{1}\boldsymbol{R}^{-1}\boldsymbol{v}_{2}=\sum_{k=1}^{\infty}\left[\boldsymbol{n}\right]_{k}^{2}\frac{\omega_{1}^{2}}{\omega_{k}^{2}}.\label{eq:StringL2inf}
\end{equation}
Note that this is $-\omega_{1}^{2}$ times the static displacement
of the string under unit load at $\xi=\xi^{\star}$. Using the expressions
for $z_{1,k}(t)$ and integrating $\left[\boldsymbol{L}(t)\right]_{2}-\left[\boldsymbol{L}^{\infty}\right]_{2}$
as per definition (\ref{eq:L1def}) we get 
\begin{multline}
\left[\boldsymbol{L}^{1}(\tau)\right]_{2}=\sum_{k=1}^{\infty}\left[\boldsymbol{n}\right]_{k}^{2}\frac{e^{-tD_{k}\omega_{k}}}{\sqrt{1-D_{k}^{2}}\omega_{k}^{3}}\left\{ 2\omega_{1}^{2}D_{k}\sqrt{1-D_{k}^{2}}\cos\left(t\sqrt{1-D_{k}^{2}}\omega_{k}\right)\right.\\
\left.+\left(\omega_{1}^{2}\left(2D_{k}^{2}-1\right)+\left(2D_{1}\omega_{1}+1\right)\omega_{k}^{2}\right)\sin\left(t\sqrt{1-D_{k}^{2}}\omega_{k}\right)\right\} -\left[\boldsymbol{n}\right]_{k}^{2}\frac{2\omega_{1}^{2}D_{k}}{\omega_{k}^{3}}.\label{eq:StringL2intGen}
\end{multline}
Assuming that $D_{k}=D$ are constant, the right limit of $\left[\boldsymbol{L}^{1}(\tau)\right]_{2}$
becomes
\begin{equation}
\lim_{\tau\to0+}\left[\boldsymbol{L}^{1}(\tau)\right]_{2}=\frac{\cos^{-1}D}{2\pi c\sqrt{1-D^{2}}}\label{eq:StringL1Lim}
\end{equation}
The detailed calculation of (\ref{eq:StringL1Lim}) can be found in
appendix \ref{sec:L1Lim}, which also indicates the boundedness of
$\left[\boldsymbol{L}^{1}(\tau)\right]_{2}$. The graph of $\left[\boldsymbol{L}^{1}(\tau)\right]_{2}$
for $c=1$ is illustrated in Fig. \ref{fig:StringL1fig}(a) for both
the conservative and the damped case.

Note that the convolution kernel $\boldsymbol{L}^{1}(\tau)$ for $D_{k}=0$
is periodic and therefore the delay that occurs as the effect of nonlinearities
is infinite. If damping is introduced $\boldsymbol{L}^{1}(\tau)$
decays in time so that the past of the system will have a smaller
effect. This is illustrated by the red line in Fig \ref{fig:StringL1fig}(a).
For non-zero damping, as an approximation one can truncate the delay
to a finite time-interval. Truncation is a reasonable choice for most
practical purposes, but it is not quite clear what are the theoretical
implications \citep{FarkasStepan}.

It is worth noting that equation (\ref{eq:StringPDE}) in the conservative
case can be solved using D'Alembert's formula that also leads to a
delay-differential equation \citep{StepanForgeDETC}, which is similar
to (\ref{eq:MZL1eq}).

\subsection{Euler-Bernoulli cantilever beam}

We choose the Euler-Bernoulli beam as our second example to illustrate
the theory. This model can support waves of infinite speed, therefore
its physical validity is questionable. Nevertheless it is worth investigating
how this property of the Euler-Bernoulli beam translates into the
properties of the memory kernel. The non-dimensional governing equation
and boundary conditions are 

\begin{equation}
\frac{\partial^{2}u}{\partial t^{2}}=-\frac{\partial^{4}u}{\partial\xi^{2}}+\mbox{damping},\; u(t,0)=u'(t,0)=u''(t,1)=u'''(t,1)=0.\label{eq:EBbeam}
\end{equation}
The natural frequencies of (\ref{eq:EBbeam}) are determined by the
equation $1+\cos\sqrt{\omega_{k}}\cosh\sqrt{\omega_{k}}=0$, which
can be approximated by $\cos\sqrt{\omega_{k}}\approx0$ for $\omega_{k}$
sufficiently large. Therefore $\omega_{k}\approx\left(k\pi-\nicefrac{\pi}{2}\right)^{2}$.
We use this estimate as a starting point to numerically find more
accurate $\omega_{k}$ values. The mode shapes at the free end of
the beam assume the values given by vector $\boldsymbol{n}=(2,-2,2,-2,\ldots)^{T}.$
We choose $\boldsymbol{m}=(\nicefrac{1}{2},0,0,\ldots)^{T}$ in formula
(\ref{eq:WmDef}), so that $\boldsymbol{W}$ satisfies our assumption
(\ref{eq:RangeCond}).

\begin{figure}
\begin{centering}
\includegraphics[width=0.49\linewidth]{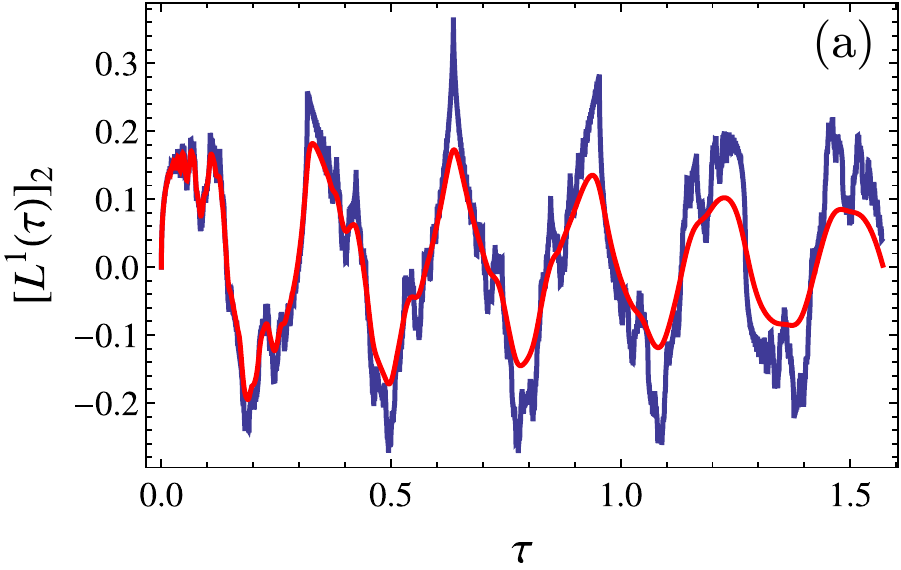}\,\includegraphics[width=0.49\linewidth]{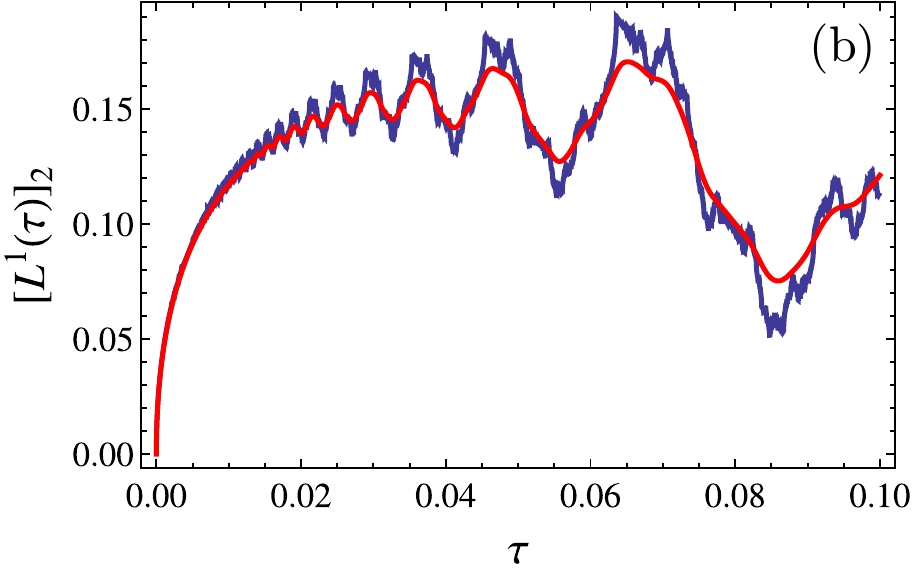}
\par\end{centering}

\caption{\label{fig:EBfig}(colour online) The graph of $\left[\boldsymbol{L}^{1}(\tau)\right]_{2}$
for the Euler-Bernoulli cantilever beam. The blue curves correspond
to the conservative case and the red curves represent the damped system
with $D_{k}=\nicefrac{1}{50}$. The conservative case illustrates
the lack of smoothness of $\left[\boldsymbol{L}^{1}(\tau)\right]_{2}$
for $\tau>0$. Panel (b) illustrates that $\left[\boldsymbol{L}^{1}(\tau)\right]_{2}$
is not differentiable at $\tau=0$. }
\end{figure}

The general formulae that were derived in section \ref{sub:StringVib}
still apply to equation (\ref{eq:EBbeam}) with the appropriate $\omega_{k}$,
$D_{k}$ and $\left[\boldsymbol{n}\right]_{k}$ values. In particular,
we use (\ref{eq:StringL2intGen}) to plot the memory kernel in Fig.
\ref{fig:EBfig}. The graph of $\left[\boldsymbol{L}^{1}(\tau)\right]_{2}$
in Fig. \ref{fig:EBfig}(a) shows that the quadratically growing natural
frequencies make the function non-smooth in the conservative case.
When damping is introduced the function becomes smooth for $\tau>0$.
Fig. \ref{fig:EBfig}(b) shows that for $0\le\tau\ll1$ $\left[\boldsymbol{L}^{1}(\tau)\right]_{2}$
grows like a power curve $\tau^{p}$, $0<p<1$. Therefore $\left[\boldsymbol{L}^{1}(\tau)\right]_{2}$
is not differentiable at $\tau=0$.

\subsection{The convergence of\textmd{ $\boldsymbol{L}_{j}^{1}(\tau)$\label{sub:L1conv}}}

So far we derived the reduced equation of motion (\ref{eq:MZL1eq})
without any consideration whether intermediate terms are well-defined.
To make the analysis rigorous we introduce infinite dimensional vector
spaces 
\begin{equation}
\boldsymbol{X}=\{\boldsymbol{x}\in\mathbb{R}^{\infty}:\sum_{i}x_{i}^{2}<\infty\}\;\mbox{and}\;\boldsymbol{Z}=\boldsymbol{X}^{2}\label{eq:Xspace}
\end{equation}
that contain the solutions of the second (\ref{eq:GenLinSystem})
and first order system (\ref{eq:1storderODE}), respectively. One
can check that in the previous two examples $\boldsymbol{v}_{M+j}\notin\boldsymbol{Z}$,
because their norm is infinite. This also implies that $\left\Vert \boldsymbol{V}\right\Vert =\infty$.
Even if we know that $\left\Vert \mathrm{e}^{\boldsymbol{R}\tau}\right\Vert \le M_{0}<\infty$,
the bound of $\boldsymbol{L}_{j}^{1}(\tau)$ cannot be directly estimated
in the straightforward way, because $\left\Vert \boldsymbol{V}\mathrm{e}^{\boldsymbol{R}\tau}\boldsymbol{v}_{M+j}\right\Vert \le\left\Vert \boldsymbol{V}\right\Vert \left\Vert \mathrm{e}^{\boldsymbol{R}\tau}\right\Vert \left\Vert \boldsymbol{v}_{M+j}\right\Vert =\infty$.
However the two examples in the previous sections show that $\boldsymbol{L}_{j}^{1}(\tau)$
can be bounded, but not necessarily smooth. In case of the string
and without damping, $\boldsymbol{L}_{j}^{1}(\tau)$ is a piecewise-smooth
function, while $\boldsymbol{L}_{j}^{1}(\tau)$ appears to be continuous
but non-differentiable for the undamped Euler-Bernoulli beam. If damping
is introduced, $\boldsymbol{L}_{j}^{1}(\tau)$ becomes smooth for
$\tau>0$ and discontinuity or non-differentiability occurs only at
$\tau=0$ in both examples.

Boundedness and smoothness depends on the eigenvalues of $\boldsymbol{R}$,
which are directly related to the natural frequencies and damping
ratios of system (\ref{eq:GenLinSystem}). First we assume that all
the eigenvalues $\sigma(\boldsymbol{R})$ are in the left half of
the complex plane including the imaginary axis, that is,
\begin{equation}
\sigma(\boldsymbol{R})\subset\{\lambda\in\mathbb{C}:\Re\lambda\le0\}.\label{eq:NegativeRealPart}
\end{equation}
This guarantees that $\mathrm{e}^{\boldsymbol{R}t}$ is a strongly
continuous semigroup with $\left\Vert \mathrm{e}^{\boldsymbol{R}t}\right\Vert \le C_{0}$
as shown in appendix \ref{sub:C0semigroup}. This does not guarantee
continuity or boundedness of $\boldsymbol{L}_{j}^{1}(\tau)$, but
it is a necessary condition. However when condition 
\begin{equation}
\left|\lambda\boldsymbol{n}_{l}\cdot(\boldsymbol{K}+\lambda\boldsymbol{C}+\lambda^{2})^{-1}\boldsymbol{n}_{j}\right|\le M_{j,l}\;\mbox{for}\;\lambda\in\{\mathbb{C}:\Re\lambda=\gamma>0\}\label{eq:ResolventCondition}
\end{equation}
is satisfied as well, appendix \ref{sub:L1bound} shows that $\boldsymbol{L}_{j}^{1}(\tau)$
is bounded.

Smoothness of $\boldsymbol{L}_{j}^{1}(\tau)$ can be guaranteed if
we replace (\ref{eq:NegativeRealPart}) and (\ref{eq:ResolventCondition})
with a stronger assumptions. We assume that there exists a $D_{0}>0$
such that the eigenvalues of $\boldsymbol{R}$ satisfy 
\begin{equation}
\sigma(\boldsymbol{R})\subset\{\lambda\in\mathbb{C}:\Re(\lambda)\le-D_{0}\left|\Im(\lambda)\right|\},\label{eq:SectorialCond}
\end{equation}
in other words the eigenvalues of $\boldsymbol{R}$ are contained
in a sector of the imaginary half-plane. Instead of (\ref{eq:ResolventCondition})
we assume that 
\begin{equation}
\left|\lambda\boldsymbol{n}_{l}\cdot(\boldsymbol{K}+\lambda\boldsymbol{C}+\lambda^{2})^{-1}\boldsymbol{n}_{j}\right|\le M_{j,l}\;\mbox{for}\;\lambda\in\{\mathbb{C}:\left|\arg\lambda\right|=\nicefrac{\pi}{2}+\delta\},\,0<\delta<\tan^{-1}D_{0}.\label{eq:SmoothCond}
\end{equation}
In case of the modal equations (\ref{eq:linModeDescr}), assumption
(\ref{eq:SectorialCond}) holds if there is a $D_{0}$ such that for
the damping ratios 
\begin{equation}
0<D_{0}\le D_{k}.\label{eq:dampingAssumption}
\end{equation}
Condition (\ref{eq:SectorialCond}) ensures that unforced vibrations
dissipate faster for higher natural frequencies which is essential
for the smoothness of solutions. According to the theory of semigroups
\citep{Pazy83} if (\ref{eq:SectorialCond}) holds the fundamental
matrix $\mathrm{e}^{\boldsymbol{R}t}$ of equation (\ref{eq:1storderODE})
is a holomorphic function of $t$, in other words, its Taylor series
converges in a sector about the non-negative real axis $t\ge0$ within
the complex plane. Appendix \ref{sub:L1bound} proves that if (\ref{eq:SectorialCond})
and (\ref{eq:SmoothCond}) are satisfied $\boldsymbol{L}_{j}^{1}(\tau)$
is smooth for $\tau>0$.

Condition (\ref{eq:ResolventCondition}) has a mechanical meaning.
When the structure is forced at contact point $\chi_{j}$ with $f_{j}(t)=\mathrm{e}^{\gamma t}\cos\omega t$,
$\gamma>0$, the velocity response $y_{M+l}(t)$ when scaled back
with the exponential growth must be bounded independent of the forcing
frequency $\omega$, that is $\left|y_{M+l}(t)e^{-\gamma t}\right|\le M_{j,l}$.
For smoothness of $\boldsymbol{L}_{j}^{1}(\tau)$ we require that
the decaying forcing $f_{j}(t)=\mathrm{e}^{-\delta\omega t}\cos\omega t$
produces a similarly decaying velocity with $\left|y_{M+l}(t)\mathrm{e}^{\delta\omega t}\right|\le M_{j,l}$
for $0<\delta<D_{0}$ independent of $\omega$.

In general, it is not straightforward to check whether (\ref{eq:ResolventCondition})
holds. Let us consider the modal equations (\ref{eq:linModeDescr})
without damping and assume that the natural frequencies scale as $\omega_{k}=\omega_{0}k^{\nicefrac{\ell}{2}},$
where $\ell=2,3,4,\dots$. We also assume that $\sup_{k}\left|\left[\boldsymbol{n}_{j}\right]_{k}\left[\boldsymbol{n}_{l}\right]_{k}\right|\le C_{j,l}$,
which implies that 
\begin{equation}
\left|\lambda\boldsymbol{n}_{l}\cdot(\boldsymbol{K}+\lambda^{2})^{-1}\boldsymbol{n}_{j}\right|\le C_{j,l}\sum_{k=1}^{\infty}\frac{\lambda}{\omega_{0}^{2}k^{\ell}+\lambda^{2}}=C_{j,l}\left(\sum_{r=1}^{\ell}\frac{\psi\left(-(-1)^{r/\ell}\lambda^{2/\ell}\right)}{\ell\omega_{0}^{2/\ell}(-1)^{r-r/\ell}\lambda^{1-2/\ell}}-\frac{1}{\lambda}\right),\label{eq:UndampResolv}
\end{equation}
where $\psi(\cdot)$ is the logarithmic derivative of the Euler Gamma
function \citep{Abramowitz}. Function $\psi$ has isolated singularities
on the negative real axis, otherwise it is bounded. Therefore a $\gamma>0$
of (\ref{eq:ResolventCondition}) can be chosen such that none of
the arguments of $\psi$ goes through these singularities. This means
that there is an $M_{j,l}$ such that $\left|\lambda\boldsymbol{n}_{l}\cdot(\boldsymbol{K}+\lambda^{2})^{-1}\boldsymbol{n}_{j}\right|\le M_{j,l}$.
In particular, for the two examples of the string and the beam we
have
\begin{equation}
\left|\lambda\boldsymbol{n}_{l}\cdot(\boldsymbol{K}+\lambda^{2})^{-1}\boldsymbol{n}_{j}\right|\le C_{j,l}\begin{cases}
\left|\frac{\pi\coth\nicefrac{\pi\lambda}{\omega_{0}}}{2\omega_{0}}-\frac{1}{2\lambda}\right|\le\frac{\pi}{2\omega_{0}}\frac{\sinh\frac{2\pi\gamma}{\omega_{0}}+1}{\cosh\frac{2\pi\gamma}{\omega_{0}}-1}+\frac{1}{2\gamma}, & \ell=2\\
\left|\frac{\sqrt[4]{-1}\pi\left(\cot\left(\frac{\sqrt[4]{-1}\pi\sqrt{\lambda}}{\sqrt{\omega_{0}}}\right)+\coth\left(\frac{\sqrt[4]{-1}\pi\sqrt{\lambda}}{\sqrt{\omega_{0}}}\right)\right)}{4\sqrt{\lambda\omega_{0}}}-\frac{1}{2\lambda}\right|, & \ell=4
\end{cases}.
\end{equation}
Note that for $\ell<2$, the sum (\ref{eq:UndampResolv}) is not uniformly
bounded, due to the $\lambda^{1-2/\ell}$ term in the denominator.

\section{Non-smooth dynamics\label{sub:NSDyn}}

We are now in the position to include the strongly nonlinear contact
forces (\ref{eq:interForce}) into the reduced model (\ref{eq:MZL1eq})
and investigate their effect. For sake of simplicity in this section
we assume a single contact force $F_{c}(\boldsymbol{y})$, so that
the governing equation becomes 
\begin{equation}
\dot{\boldsymbol{y}}(t)=\boldsymbol{A}\boldsymbol{y}(t)+\boldsymbol{L}^{\infty}F_{c}(\boldsymbol{y}(t))+\int_{0}^{t-s}\mathrm{d}_{\tau}\boldsymbol{L}^{1}(\tau)\frac{\mathrm{d}}{\mathrm{d}t}[F_{c}(\boldsymbol{y}(t-\tau))]+\boldsymbol{g}(t),\label{eq:MZNSeq}
\end{equation}
where $\boldsymbol{g}(t)=\boldsymbol{H}(t)\boldsymbol{z}(s)+\boldsymbol{L}^{0}(t-s)F_{c}(\boldsymbol{y}(s)).$
The properties of solutions of (\ref{eq:MZNSeq}) strongly depend
on both $\boldsymbol{L}^{1}(\tau)$ and $F_{c}(\boldsymbol{y})$.
We also assume that $\boldsymbol{L}^{1}(\tau)$ is smooth for $\tau>0$
as it is outlined in section \ref{sub:L1conv}.

Our definition of solution at the discontinuities of $F_{c}(\boldsymbol{y})$
is based on a mechanical analogy. If the elastic bodies stick together
there is an algebraic constraint that restricts the trajectories to
sticking motion and one must be able to calculate the contact force
implicitly from equation (\ref{eq:MZNSeq}). If these contact forces
are admissible by physical law, the bodies will stick, otherwise they
will continue slipping.

To formalise this definition, we assume that $F_{c}(\boldsymbol{y})$
is discontinuous along a smooth surface defined by $h(\boldsymbol{y})=0$,
which stands for the algebraic constraint of sticking. We call $\Sigma=\{\boldsymbol{y}\in\mathbb{R}^{2M}:h(\boldsymbol{y})=0\}$
the switching surface. The physical bound of the contact force can
be defined as the two limits of $F_{c}(\boldsymbol{y})$ on the two
sides of $\Sigma$, that is,
\begin{equation}
\forall\boldsymbol{y}\in\Sigma,\quad F_{c}^{+}(\boldsymbol{y})=\lim_{\bar{\boldsymbol{y}}\to\boldsymbol{y},h(\bar{\boldsymbol{y}})>0}F_{c}(\bar{\boldsymbol{y}})\;\mbox{and}\; F_{c}^{-}(\boldsymbol{y})=\lim_{\bar{\boldsymbol{y}}\to\boldsymbol{y},h(\bar{\boldsymbol{y}})<0}F_{c}(\bar{\boldsymbol{y}}).\label{eq:FcLimits}
\end{equation}
Without restricting generality we assume that $F_{c}^{-}(\boldsymbol{y})<F_{c}^{+}(\boldsymbol{y})$.
Alternatively, $F_{c}^{-}$ and $F_{c}^{+}$ can be defined on the
switching surface $\Sigma$ independently of $F_{c}$, when one wants
to distinguish between static and dynamic friction.

According to our physical interpretation of the solution, when a trajectory
reaches the switching surface $\Sigma$ the trajectory either crosses
$\Sigma$ or becomes part of $\Sigma$, which means sticking in the
physical sense. The algebraic constraint of sticking is $h(\boldsymbol{y}(t))=0$.
While sticking the contact force $F_{c}^{\star}(t)$ must stay within
physical bounds 
\begin{equation}
F_{c}^{-}(\boldsymbol{y})\le F_{c}^{\star}(t)\le F_{c}^{+}(\boldsymbol{y})\label{eq:FcBound}
\end{equation}
and the vector field must be tangential to $\Sigma$, that is, 
\begin{equation}
\nabla h(\boldsymbol{y}(t))\cdot\dot{\boldsymbol{y}}(t)=0,\label{eq:NsSlideCond}
\end{equation}
so that the solution continues on the switching surface. If such a
contact force cannot be found the solution crosses the switching surface
and a discontinuity develops in the contact force. To calculate the
contact force $F_{c}^{\star}(t)$ that makes the solution restricted
to $\Sigma$ we substitute (\ref{eq:MZNSeq}) into (\ref{eq:NsSlideCond}),
which yields 
\begin{equation}
0=\nabla h(\boldsymbol{y}(t))\cdot\biggl\{\boldsymbol{A}\boldsymbol{y}(t)+\boldsymbol{L}^{\infty}F_{c}^{\star}(t)+\int_{0}^{t-s}\mathrm{d}_{\tau}\boldsymbol{L}^{1}(\tau)\frac{\mathrm{d}}{\mathrm{d}t}[F_{c}^{\star}(t-\tau)]+\boldsymbol{g}(t)\biggr\}.\label{eq:LambdaEq}
\end{equation}
Equation (\ref{eq:LambdaEq}) involves the history of the contact
force, which is either $F_{c}^{\star}(t)=F_{c}(\boldsymbol{y}(t))$
if $h(\boldsymbol{y}(t))\neq0$ or it is calculated from (\ref{eq:LambdaEq}).

The question is whether the contact force $F_{c}^{\star}(t)$ is well
defined during the stick phase by equation (\ref{eq:LambdaEq}), which
is an integral equation for $F_{c}^{\star}(t)$. To answer this we
need to consider possible singularities of $\boldsymbol{L}^{1}(\tau)$
at $\tau=0$. Since $\boldsymbol{L}^{1}(\tau)$ is bounded one can
find a maximal $0\le\alpha\le1$ and a positive constant $C$ such
that $\left\Vert \boldsymbol{L}^{1}(\tau)\right\Vert <C\tau^{\alpha}$.
This is called the H\"older condition and $\alpha$ is the H\"older
exponent. If $\alpha<1$ we can also find a constant $\boldsymbol{L}^{1+}$
and positive $C_{0}$ such that 
\[
\left\Vert \boldsymbol{L}^{1}(\tau)-\boldsymbol{L}^{1+}\tau^{\alpha}\right\Vert <C_{0}\tau.
\]
This means that $\boldsymbol{L}^{1}(\tau)$ is a sum of the singular
$\boldsymbol{L}^{1+}\tau^{\alpha}$ and a differentiable function.
There are three cases to consider:
\begin{enumerate}
\item $\alpha=1$, so that $\boldsymbol{L}^{1}(\tau)$ is differentiable.
We assume that $\nabla h(\boldsymbol{y}(t))\cdot\boldsymbol{L}(0)\neq0$.\label{enu:C0}
\item $\alpha=0$, so that $\boldsymbol{L}^{1}(\tau)$ is discontinuous
and $\boldsymbol{L}^{1+}=\lim_{\tau\to0+}\boldsymbol{L}^{1}(\tau)$.
We assume that $\nabla h(\boldsymbol{y}(t))\cdot\boldsymbol{L}^{1+}\neq0$.\label{enu:C1}
\item $0<\alpha<1$, when $\boldsymbol{L}^{1}(\tau)$ is not differentiable,
but continuous. Similarly, we assume that $\nabla h(\boldsymbol{y}(t))\cdot\boldsymbol{L}^{1+}\neq0$.\label{enu:C2}
\end{enumerate}
In case \ref{enu:C0}, when $\boldsymbol{L}^{1}(\tau)$ is continuously
differentiable on $[0,\infty)$, the integral term can be expressed
using $\boldsymbol{L}(\tau)$ as in equation (\ref{eq:MZequation}).
This is the case when the governing equations are finite dimensional
or $\boldsymbol{n}_{j}$ have finite norms. Therefore the same dynamical
phenomena should occur as in finite dimensional systems, which cannot
be resolved by our method. Applying (\ref{eq:NsSlideCond}) to equation
(\ref{eq:MZequation}) we find that the contact force obeys the integral
equation 

\begin{equation}
F_{c}^{\star}(t)=\frac{-\nabla h(\boldsymbol{y}(t))}{\nabla h(\boldsymbol{y}(t))\cdot\boldsymbol{L}(0)}\cdot\biggl\{\boldsymbol{A}\boldsymbol{y}(t)+\int_{0}^{t-s}\mathrm{d}_{\tau}\boldsymbol{L}(\tau)F_{c}^{\star}(t-\tau)+\boldsymbol{H}(t)\boldsymbol{z}(s)\biggr\}.\label{eq:FcL1cont}
\end{equation}
Due to the differentiability of $\boldsymbol{L}^{1}(\tau)$ its derivative
$\boldsymbol{L}^{0}(\tau)$ and therefore $\boldsymbol{L}(0)$ must
be bounded. When calculating the contact force by equation (\ref{eq:FcL1cont})
can result in a discontinuity of $F_{c}^{\star}(t)$ at the onset
of the stick phase. Another cause of singularity is when $\nabla h(\boldsymbol{y}(t))\cdot\boldsymbol{L}(0)=0$,
which can occur in case of the two-fold singularity \citep{ColomboJeffrey}.

The pre-tensed string model falls into case \ref{enu:C1}. Due to
the discontinuity of $\boldsymbol{L}^{1}(\tau)$, equation (\ref{eq:LambdaEq})
can be rearranged as a delay differential equation with $\nicefrac{\mathrm{d}}{\mathrm{d}t}F_{c}^{\star}(t)$
on the left-hand side, that is, 
\begin{equation}
\frac{\mathrm{d}}{\mathrm{d}t}F_{c}^{\star}(t)=\frac{-\nabla h(\boldsymbol{y}(t))}{\nabla h(\boldsymbol{y}(t))\cdot\boldsymbol{L}^{1+}}\cdot\biggl\{\boldsymbol{A}\boldsymbol{y}(t)+\boldsymbol{L}^{\infty}F_{c}^{\star}(t)+\int_{0+}^{t-s}\mathrm{d}_{\tau}\boldsymbol{L}^{1}(\tau)\frac{\mathrm{d}}{\mathrm{d}t}[F_{c}^{\star}(t-\tau)]+\boldsymbol{g}(t)\biggr\}.\label{eq:FcL1Jump}
\end{equation}
At the onset of stick at $t^{\star}$ the initial condition is $F_{c}^{\star}(t^{\star})=\lim_{t\to t^{\star}-0}F_{c}(\boldsymbol{y}(t))$.
Since all the terms in (\ref{eq:FcL1Jump}) are bounded $\nicefrac{\mathrm{d}}{\mathrm{d}t}F_{c}^{\star}(t^{\star})$
must also be bounded. Therefore $F_{c}^{\star}(t)$ is a Lipschitz
continuous function of time when the solution gets restricted to $\Sigma$
and all throughout the stick phase. At the transition from stick to
slip $F_{c}^{\star}(t)$ is continuous if $F_{c}^{\pm}(\boldsymbol{y})$
is the limit of $F_{c}(\boldsymbol{y})$ defined by (\ref{eq:FcLimits}).
If in addition the slope of $F_{c}(\boldsymbol{y})$ on the relevant
side of $\Sigma$ is finite, $F_{c}^{\star}(t)$ is Lipschitz continuous.
It remains to be investigated what are the dynamical consequences
when $\nabla h(\boldsymbol{y}(t))\cdot\boldsymbol{L}^{1+}=0$ and
whether the uniqueness of solution is preserved through such a singularity.

The Euler-Bernoulli beam falls into case \ref{enu:C2}. First we note
that 
\begin{equation}
\int_{0}^{t-s}\mathrm{d}_{\tau}\tau^{\alpha}\frac{\mathrm{d}}{\mathrm{d}t}[F_{c}^{\star}(t-\tau)]=\int_{s}^{t}(t-\theta)^{\alpha-1}\frac{\mathrm{d}}{\mathrm{d}\theta}[F_{c}^{\star}(\theta)]\mathrm{d}\theta,
\end{equation}
which is by definition $-\Gamma(\alpha)$ times the $\alpha$ fractional
integral of $\nicefrac{\mathrm{d}}{\mathrm{d}t}F_{c}^{\star}(t)$
\citep{mcbride1979fractional}. We assume that the stick phase starts
at time $t^{\star}$. Using the rules of fractional integration we
find that 
\begin{equation}
\int_{t^{\star}}^{t}(t-\theta)^{-\alpha}\left(\int_{0}^{\theta-t^{\star}}\mathrm{d}_{\tau}\tau^{\alpha}\frac{\mathrm{d}}{\mathrm{d}t}[F_{c}^{\star}(\theta-\tau)]\right)\mathrm{d}\theta=\frac{\alpha\pi}{\sin\alpha\pi}\left(F_{c}^{\star}(t)-F_{c}^{\star}(t^{\star})\right).\label{eq:FracIntTest}
\end{equation}
By separating the singular component of equation (\ref{eq:LambdaEq})
we get
\begin{multline}
\int_{0}^{t-t^{\star}}\mathrm{d}_{\tau}\tau^{\alpha}\frac{\mathrm{d}}{\mathrm{d}t}[F_{c}^{\star}(t-\tau)]=\frac{-\nabla h(\boldsymbol{y}(t))}{\nabla h(\boldsymbol{y}(t))\cdot\boldsymbol{L}^{1+}}\cdot\biggl(\boldsymbol{A}\boldsymbol{y}(t)+\boldsymbol{L}^{\infty}F_{c}^{\star}(t)+\\
+\int_{0}^{t-t\star}\mathrm{d}_{\tau}\left(\boldsymbol{L}^{1}(\tau)-\boldsymbol{L}^{1+}\tau^{\alpha}\right)\frac{\mathrm{d}}{\mathrm{d}\theta}[F_{c}^{\star}(t-\tau)]+\int_{t-t^{\star}}^{t-s}\mathrm{d}_{\tau}\boldsymbol{L}^{1}(\tau)\frac{\mathrm{d}}{\mathrm{d}t}[F_{c}^{\star}(t-\tau)]+\boldsymbol{g}(t)\biggr).\label{eq:SingSep}
\end{multline}
Since we assumed that $\nabla h(\boldsymbol{y}(t))\cdot\boldsymbol{L}^{1+}\neq0$,
it follows that all the terms on the right side of (\ref{eq:SingSep})
are bounded by $C_{1}$. Therefore we fractional integrate (\ref{eq:SingSep})
with $1-\alpha$ exponent exactly as in (\ref{eq:FracIntTest}) and
get
\begin{equation}
\left|\frac{\alpha\pi}{\sin\alpha\pi}\left(F_{c}^{\star}(t)-F_{c}^{\star}(t^{\star})\right)\right|\le\int_{t^{\star}}^{t}(t-\theta)^{-\alpha}C_{1}\mathrm{d}\theta=\frac{C_{1}\left(t-t^{\star}\right)^{1-\alpha}}{1-\alpha}.\label{eq:FcHolder}
\end{equation}
This means that if $\nabla h(\boldsymbol{y}(t))\cdot\boldsymbol{L}^{1+}\neq0$,
$F_{c}^{\star}(t)$ is H\"older continuous with exponent $1-\alpha$.
We note that H\"older continuity implies continuity in the traditional
sense, therefore the friction force is continuos during the transition
from slip to stick.

\section{Stick-slip motion of a bowed string\label{sec:Bowstring}}

\begin{figure}
\begin{centering}
\includegraphics[width=0.9\linewidth]{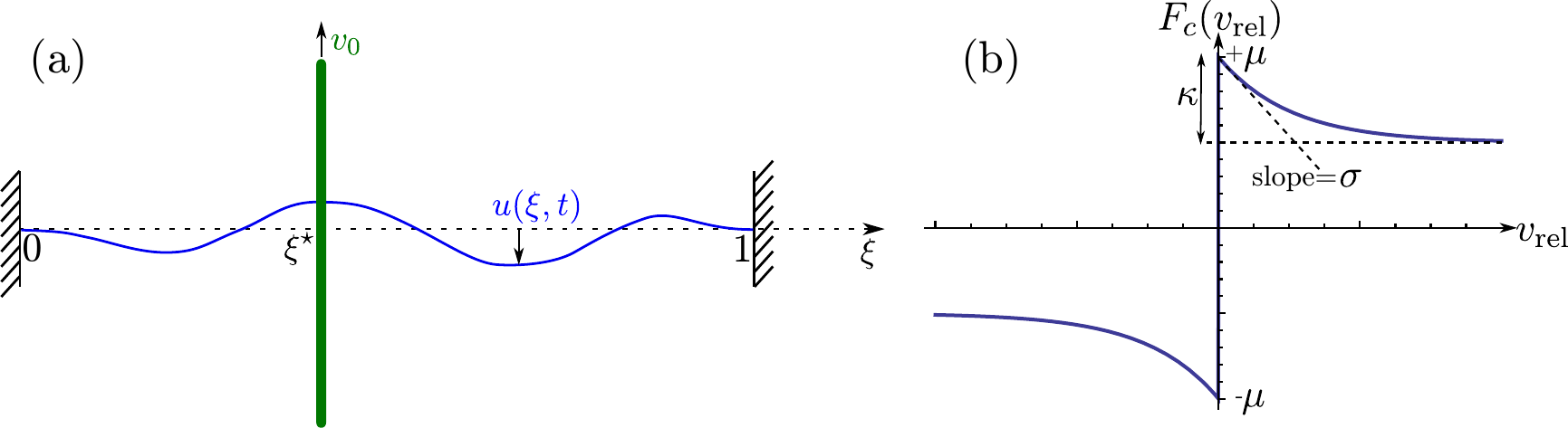}
\par\end{centering}

\caption{(colour online) (a) Schematic of a bowed string. The bow is pulled
with a constant velocity $v_{0}$, while the string exhibits a stick-slip
vibration generated by the friction between the bow and the string.\label{fig:bowedstring}
(b) Graph of the Coulomb-like friction force.}
\end{figure}
To see our theory applied to a mechanical system consider the example
of a bowed string in Fig. \ref{fig:bowedstring}(a). We consider the
same equation of motion as in section \ref{sub:StringVib}, where
we derived all the necessary ingredients of the reduced model apart
from the contact force. Because $\boldsymbol{L}^{1}(\tau)$ has a
discontinuity at $\tau=0$ this example falls into case \ref{enu:C1}.~of
section \ref{sub:NSDyn}.

To complete the model we define the contact force $F_{c}$ of equation
(\ref{eq:MZNSeq}) as the friction force between the bow and the string.
We assume that the string is being bowed at $\xi=\xi^{\star}$ with
velocity $v_{0}$ that generates the friction force 
\begin{equation}
f_{c}(v_{\textrm{rel}})=\mathrm{sign}v_{\mathrm{rel}}\left(\mu-\kappa+\kappa\exp\left(-\sigma\left|v_{\mathrm{rel}}\right|\right)\right),\label{eq:fforce}
\end{equation}
where $v_{\textrm{rel}}$ is the relative velocity of the string and
the bow. The graph of the friction force function can be seen in Fig.
\ref{fig:bowedstring}(b). In this example the static friction force
is within the interval $[-\mu,\mu]$. The relative velocity between
the string and the bow is expressed as a function of a resolved variable
$v_{\mathrm{rel}}=y_{2}(t)-v_{0}$. We also use the relative velocity
to define the switching surface $\Sigma$ by $h(\boldsymbol{y}(t))=v_{\mathrm{rel}}$.
Therefore the contact force of equation (\ref{eq:MZNSeq}) becomes
$F_{c}(\boldsymbol{y}(t))=f_{c}(h(\boldsymbol{y}(t)))$.

\subsection{Numerical method\label{sub:NumMeth}}

We use a simple explicit Euler method to approximate the solutions
of (\ref{eq:MZNSeq}) and (\ref{eq:LambdaEq}). We assume that time
is quantised in $\varepsilon$ chunks, so that $\boldsymbol{y}_{q}=\boldsymbol{y}(q\varepsilon)$,
$f_{c,q}=f_{c}(q\varepsilon)$, where $q=0,1,2,\ldots$. In case of
slipping the only unknown is the state variable $\boldsymbol{y}_{q}$
that is calculated using the formula
\begin{multline}
\boldsymbol{y}_{q+1}=\boldsymbol{y}_{q}+\varepsilon\Biggl(\boldsymbol{A}\boldsymbol{y}_{q}+\boldsymbol{L}^{\infty}f_{c,q}+\sum_{j=0}^{n}\boldsymbol{L}^{0}(j\varepsilon)\left(f_{c,q-j}-f_{c,q-j-1}\right)\\
+\boldsymbol{L}^{0}(q\varepsilon)f_{c,0}+\boldsymbol{H}(q\varepsilon)\boldsymbol{z}(0)\Biggr),\label{eq:NumStateUpdate}
\end{multline}
where the friction force $f_{c,q}=F_{c}\left(\left[\boldsymbol{y}_{q}\right]_{2}-v_{0}\right)$
is used. The integration is approximated by the rectangle rule. For
just illustrating the theory such a crude approximation is sufficient
while for better accuracy and efficiency higher order methods, such
as the Runge-Kutta \citep{Iserles} method could be used. In our calculations
we keep the step size reasonably short at $\varepsilon=5\times10^{-4}$.

If the relative velocity $h(\boldsymbol{y}_{q})=\left[\boldsymbol{y}_{q}\right]_{2}-v_{0}$
of the string and the bow becomes zero there are two possibilities.
Either the trajectory crosses the switching surface $\Sigma$ or it
will stay on $\Sigma$ satisfying the equation 
\begin{equation}
\nabla h(\boldsymbol{y}_{q})\cdot(\boldsymbol{y}_{q+1}-\boldsymbol{y}_{q})=0,\label{eq:NumTangency}
\end{equation}
that is, the discretised version of $h(\boldsymbol{y}(t))\cdot\dot{\boldsymbol{y}}(t)=0$.
To test which case applies, we substitute (\ref{eq:NumStateUpdate})
into (\ref{eq:NumTangency}) and solve for the friction force $f_{c,q}$
that would hold the string and the bow together, which becomes 
\begin{multline}
f_{c,q}=\frac{-\nabla h(\boldsymbol{y}_{q})}{\nabla h(\boldsymbol{y}_{q})\cdot(\boldsymbol{L}^{\infty}+\boldsymbol{L}^{0}(0))}\Biggl(-\boldsymbol{L}^{0}(0)f_{c,q-1}+\boldsymbol{A}\boldsymbol{y}_{q}+\sum_{j=1}^{q}\boldsymbol{L}^{0}(j)\left(f_{c,q-j}-f_{c,q-j-1}\right)\\
+\boldsymbol{L}^{0}(q\varepsilon)f_{c,0}+\boldsymbol{H}(q\varepsilon)\boldsymbol{z}(0)\Biggr).\label{eq:NumFcUpdate}
\end{multline}
If the calculated friction force satisfies $-\mu\le f_{c,q}\le\mu$,
the bow and the string stick together. For the stick phase of motion
we use equation (\ref{eq:NumStateUpdate}) to advance the solution
together with this dynamic friction force of equation (\ref{eq:NumFcUpdate}).

\subsection{Numerical results}

\begin{figure}
\begin{centering}
\includegraphics[width=0.99\linewidth]{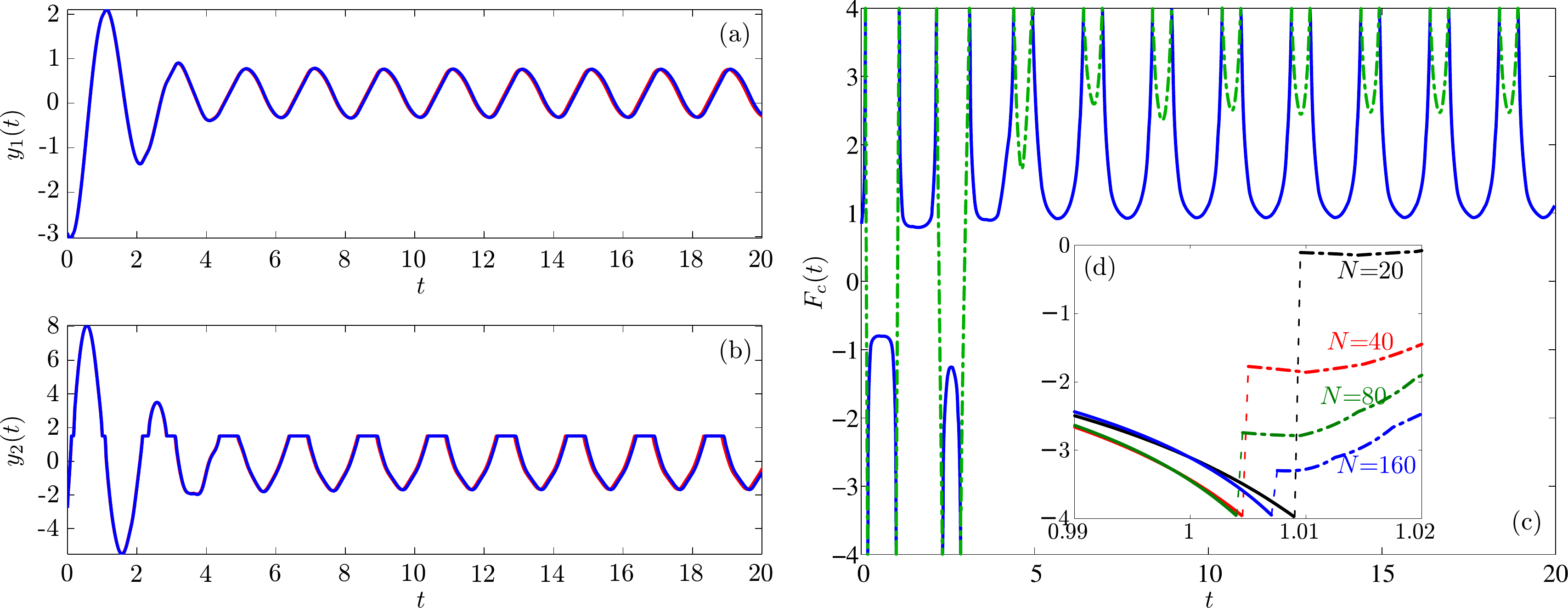}
\par\end{centering}

\caption{Solution trajectories of equations (\ref{eq:linModeDescr}) and (\ref{eq:MZNSeq}).
(a) Displacement of the string and (b) velocity of the string at the
contact point. (c) Friction force between the bow and string using
equation (\ref{eq:MZNSeq}). The continuous lines denote slip, the
dash-dotted lines represent sticking motion. (d) The discontinuity
of the friction force at the onset of sticking disappears in the continuum
limit. $N$ indicates the considered number of modes.\label{fig:Solution}}
\end{figure}

To illustrate the properties of the dimension reduced equation (\ref{eq:MZNSeq})
we calculated a typical stick-slip trajectory starting at a the initial
condition $\boldsymbol{z}(0)=y_{1}(0)\boldsymbol{w}_{1}+y_{2}(0)\boldsymbol{w}_{2}$
with $y_{1}(0)=-2.9224$ and $y_{2}(0)=-2.7668$. The parameters of
the friction force in equation (\ref{eq:fforce}) are $\mu=4$, $\kappa=0.32$,
$\sigma=1$, the speed of the bow is $v_{0}=\nicefrac{3}{2}$, the
damping ratios are $D_{k}=\nicefrac{1}{10}$, $k=1,\ldots,N$ and
the wave speed on the string is $c=1$. The string is bowed at $\xi^{\star}=0.4$.%
\footnote{The choice of these parameters was guided by the desire of producing
stick-slip motion rather than physical consideration.%
} We solved equation (\ref{eq:linModeDescr}) using \noun{Matlab}'s
\textsf{ode113} solver and (\ref{eq:MZNSeq}) using our method described
in section \ref{sub:NumMeth}. The results of the the simulation for
reduced and the full model shown in Figure \ref{fig:Solution}(a,b)
are nearly indistinguishable, because the only approximations are
within the numerical methods. The blue curves denote the solution
of (\ref{eq:MZNSeq}) and the (almost invisible) red curve underneath
represents the solution of (\ref{eq:linModeDescr}) using solution
techniques described by \citet{PiiroinenSim}.

Initially the solution spends short time intervals on the switching
surface and then settles into a periodic stick-slip motion. The stick
phases can be recognised in Fig. \ref{fig:Solution}(b) as short horizontal
sections at $y_{2}=1.5$. In Fig. \ref{fig:Solution}(c) the friction
force is represented by the blue lines and the green dash-dotted lines
for the slipping and the sticking motion, respectively. The friction
force also appears to be discontinuous. To calculate this solution
we did not use the converged $\boldsymbol{L}^{1}(\tau)$, instead
we used a series of mode truncations shown in Fig.~(\ref{fig:StringL1fig})(b).
On a smaller scale Fig. \ref{fig:Solution}(d) shows that the gap
in the friction force (dashed line) between the slipping segment (continuous
line) and the sticking segment (dash dotted line) of the friction
force vanishes as increasing number of modes of system (\ref{eq:linModeDescr})
are considered. As the theory dictates the gap should vanish in the
infinite dimensional case.

\section{Conclusions}

In this paper we considered vibrations of structures that are composed
of linear elastic bodies coupled through strongly nonlinear contact
forces such as friction. The coupling was assumed to occur at point
contacts. We introduced an exact transformation based on the Mori-Zwanzig
formalism that reduces the infinite dimensional system of ordinary
differential equations to a description with time delay involving
small number of variables. We found that the model reduction technique
converges and contact forces become continuous even though the governing
equation is discontinuous. We illustrated this novel technique through
the example of a bowed string.

Through examples we found that if natural frequencies scale linearly
with the mode number, the contact forces are Lipschitz continuous
during the transition from slip to stick. This is the case of the
elastic string. If the natural frequencies increase faster than linear,
the contact forces are only continuous. The Euler-Bernoulli beam exhibits
such a behaviour, but it also allows infinite wave speed, which can
be though of as not physical. In reality however, every structure
will have small scale longitudinal vibration components with linearly
scaled frequencies similar to the Timoshenko beam model \citep{SzalaiMZImpact}.
We expect that if all the details are considered for a linear structure,
the contact forces must always be Lipschitz continuous in time. This
finding together with the new form of governing equations could be
used in further studies to understand the source of non-deterministic
motion \citep{ColomboJeffrey}.

The reduced equations are also structurally stable. Small perturbations
to the memory kernel or other terms only deform solutions but do not
change their qualitative behaviour as long as the qualitative features
of the memory kernel are preserved. This is a clear advantage over
finite dimensional approximation of non-smooth systems, where small
perturbations can cause qualitatively different solutions. Therefore
once the qualitative form of the memory kernel is established non-smooth
mechanical systems can be approximated more successfully using our
description.

How non-smooth phenomena of low dimensional systems manifest in continuum
structures is an open question. For low dimensional systems many singularities
can occur that lead to chaotic and resonant vibration on invariant
polygons \citep{SzalaiPolygons}, the Painleve paradox \citep{PainleveAlan}
and other types of discontinuity induced bifurcations. It remains
to investigate how these phenomena occur in systems involving elastic
structures and hence equations with memory.

Our theory is developed for linear structures coupled to strong nonlinearities.
It is however possible to extend this framework to cases where the
underlying structure is nonlinear. For the weakly nonlinear case the
Hartman-Grobman theorem \citep{CoddingtonLevinson,KuznetsovBook}
guarantees the existence of a transformation that takes any weakly
nonlinear system into a linear system about an equilibrium if that
system is not undergoing a stability change. This generalisation is
currently being worked on by the author.

We also assumed point contacts in our derivations. This is a significant
simplification since most contact problems occur along a surface.
The difficulty arises when one needs to deal with contacting surfaces
that slip at one part of the contact surface while stick at others.
An interesting question is if it is possible to develop a similar
model reduction technique of such problems to involve only finite
number of variables.

\section*{Acknowledgements}

The author would like to thank G{\'a}bor St{\'e}p{\'a}n, who brought
his attention to the work of \citet{ChorinPNAS}. He would also like
thank Jan Sieber, Alan R. Champneys and John Hogan for useful discussion
and comments on the manuscript.

\appendix

\section{Model transformation\label{sec:AppTrafo}}

In this appendix we show that the infinite dimensional system (\ref{eq:1storderODE})
can be transformed into a finite dimensional delayed equation. The
delay equation involves convolution integrals that can be related
to Green's functions, but only for part of the system. The procedure
is based on the variation-of-parameters formula.

Consider the following linear forced system

\begin{equation}
\dot{\boldsymbol{z}}(t)=\boldsymbol{R}\boldsymbol{z}(t)+\boldsymbol{v}f(t),\label{eq:AppLinEq}
\end{equation}
where $\boldsymbol{z}(t),\boldsymbol{v}\in\mathbb{R}^{2N}$, $f(t)\in\mathbb{R}$
and $\boldsymbol{R}\in\mathbb{R}^{2N\times2N}$. Assume matrices $\boldsymbol{V}\in\mathbb{R}^{2M\times2N}$
and $\boldsymbol{W}\in\mathbb{R}^{2N\times2M}$ such that $\boldsymbol{S}=\boldsymbol{W}\boldsymbol{V}$
is a projection matrix $\boldsymbol{S}=\boldsymbol{S}^{2}$ with a
$2M$ dimensional range. $\boldsymbol{S}$ is a projection if and
only if $\boldsymbol{V}\boldsymbol{W}=\boldsymbol{I}_{2M}$, where
$\boldsymbol{I}_{2M}$ is the $2M$ dimensional identity. Also, define
the complementary projection matrix $\boldsymbol{Q}=\boldsymbol{I}-\boldsymbol{S}$
and the resolved coordinates $\boldsymbol{y}(t)=\boldsymbol{V}\boldsymbol{z}(t)$.
With this notation we rewrite equation (\ref{eq:AppLinEq}) into 
\begin{equation}
\dot{\boldsymbol{z}}(t)=\boldsymbol{R}\boldsymbol{S}\boldsymbol{z}(t)+\boldsymbol{R}\boldsymbol{Q}\boldsymbol{z}(t)+\boldsymbol{v}f(t).\label{eq:AppSepEq}
\end{equation}
Assume that the solution of $\dot{\boldsymbol{z}}(t)=\boldsymbol{R}\boldsymbol{Q}\boldsymbol{z}(t)$
can be computed for specific initial conditions. Therefore the solution
of (\ref{eq:AppSepEq}) can formally be expressed using the variation-of-parameters
or Dyson's \citep{CoddingtonLevinson} formula as
\begin{equation}
\boldsymbol{z}(t)=\mathrm{e}^{\boldsymbol{R}\boldsymbol{Q}(t-s)}\boldsymbol{z}(s)+\int_{0}^{t-s}\mathrm{e}^{\boldsymbol{R}\boldsymbol{Q}\tau}\left(\boldsymbol{R}\boldsymbol{S}\boldsymbol{z}(t-\tau)+\boldsymbol{v}f(t-\tau)\right)\mathrm{d}\tau.
\end{equation}
Substituting this result into the second term on the right side of
(\ref{eq:AppSepEq}) we get
\begin{multline}
\dot{\boldsymbol{z}}(t)=\boldsymbol{R}\boldsymbol{S}\boldsymbol{z}(t)+\boldsymbol{v}f(t)\\
+\boldsymbol{R}\boldsymbol{Q}\left\{ \mathrm{e}^{\boldsymbol{R}\boldsymbol{Q}(t-s)}\boldsymbol{z}(s)+\int_{0}^{t-s}\mathrm{e}^{\boldsymbol{R}\boldsymbol{Q}\tau}\left(\boldsymbol{R}\boldsymbol{S}\boldsymbol{z}(t-\tau)+\boldsymbol{v}f(t-\tau)\right)\right\} \mathrm{d}\tau.\label{eq:AppSubs}
\end{multline}
Note that $\boldsymbol{S}\boldsymbol{z}(t)=\boldsymbol{W}\boldsymbol{V}\boldsymbol{z}(t)=\boldsymbol{W}\boldsymbol{y}(t)$
and project (\ref{eq:AppSubs}) using $\boldsymbol{V}$, to get 
\begin{multline}
\dot{\boldsymbol{y}}(t)=\boldsymbol{V}\boldsymbol{R}\boldsymbol{W}\boldsymbol{y}(t)+\boldsymbol{V}\boldsymbol{v}f(t)\\
+\boldsymbol{V}\boldsymbol{R}\boldsymbol{Q}\left\{ \mathrm{e}^{\boldsymbol{R}\boldsymbol{Q}(t-s)}\boldsymbol{z}(s)+\int_{0}^{t-s}\mathrm{e}^{\boldsymbol{R}\boldsymbol{Q}\tau}\left(\boldsymbol{R}\boldsymbol{S}\boldsymbol{z}(t-\tau)+\boldsymbol{v}f(t-\tau)\right)\right\} \mathrm{d}\tau,
\end{multline}
which is the reduced equation for only the resolved coordinates $\boldsymbol{y}(t)$.
Note that $\boldsymbol{R}\boldsymbol{Q}\mathrm{e}^{\boldsymbol{R}\boldsymbol{Q}\tau}=\nicefrac{d}{d\tau}\mathrm{e}^{\boldsymbol{R}\boldsymbol{Q}\tau}$
hence the integrals can be rewritten with Riemann-Stieltjes integrals
as
\begin{equation}
\dot{\boldsymbol{y}}(t)=\boldsymbol{A}\boldsymbol{y}(t)+\int_{0}^{t-s}\mathrm{d}_{\tau}\boldsymbol{K}(\tau)\boldsymbol{y}(t-\tau)+\int_{0}^{t-s}\mathrm{d}_{\tau}\boldsymbol{L}(\tau)f(t-\tau)+\boldsymbol{V}\boldsymbol{v}f(t)+\boldsymbol{H}(t)\boldsymbol{z}(s),\label{eq:FullMZeq}
\end{equation}
where
\begin{align}
\boldsymbol{A} & =\boldsymbol{V}\boldsymbol{R}\boldsymbol{W}, &  & \in\mathbb{R}^{2M\times2M}\label{eq:A-Aexpr}\\
\boldsymbol{H}(t)\boldsymbol{z}(s) & =\boldsymbol{V}\boldsymbol{R}\boldsymbol{Q}\mathrm{e}^{\boldsymbol{R}\boldsymbol{Q}(t-s)}\boldsymbol{z}(s), &  & \in\mathbb{R}^{2M}\label{eq:A-Hexpr}\\
\boldsymbol{K}(\tau) & =\boldsymbol{V}\mathrm{e}^{\boldsymbol{R}\boldsymbol{Q}\tau}\boldsymbol{R}\boldsymbol{W}, &  & \in\mathbb{R}^{2M\times2M}\label{eq:A-Kexpr}\\
\boldsymbol{L}(\tau) & =\boldsymbol{V}\mathrm{e}^{\boldsymbol{R}\boldsymbol{Q}\tau}\boldsymbol{v}. &  & \in\mathbb{R}^{2M}\label{eq:A-Lexpr}
\end{align}
If the range of $\boldsymbol{W}$ is invariant under $\boldsymbol{R}$,
then $\boldsymbol{K}(\tau)=\boldsymbol{A}$ constant. This occurs
because the image of the range of $\boldsymbol{W}$ is in the kernel
of $\boldsymbol{Q}$. Consequently the integral with $\boldsymbol{K}(\tau)$
vanishes.

Note that this procedure is a simplified version of the Mori-Zwanzig
formalism \citep{ChorinPNAS,EvansMorriss} for linear systems. Therefore
our procedure can be extended to nonlinear systems. In the nonlinear
case $\boldsymbol{A}\boldsymbol{y}$, $\boldsymbol{K}(\tau)\boldsymbol{y}$
and $\boldsymbol{L}(\tau)f$ become nonlinear functions of $\boldsymbol{y}$
and $f$, respectively.

\section{The memory kernels\label{sec:AppMemKer}}

In this appendix we show that the memory kernel can be obtained from
the solution of the first order system (\ref{eq:1storderODE}) if
condition (\ref{eq:RangeCond}) holds. Condition (\ref{eq:RangeCond})
implies that there is a $2M\times2M$ matrix $\boldsymbol{X}$ such
that $\boldsymbol{R}\boldsymbol{W}=\boldsymbol{W}\boldsymbol{X}$.
If we multiply this expression by $\boldsymbol{V}$ from the left
we get the identity $\boldsymbol{A}=\boldsymbol{V}\boldsymbol{R}\boldsymbol{W}=\boldsymbol{V}\boldsymbol{W}\boldsymbol{X}=\boldsymbol{X}$.
As a consequence for any integer $p$ we must have
\begin{equation}
\boldsymbol{R}^{p}\boldsymbol{W}=\boldsymbol{W}\boldsymbol{A}^{p}.\label{eq:WApower}
\end{equation}

To investigate the expression $\mathrm{e}^{\boldsymbol{R}\boldsymbol{Q}t}$
that occurs in the definition of the memory kernel (\ref{eq:A-Lexpr})
we define
\begin{equation}
\boldsymbol{\varPhi}(t)=\mathrm{e}^{\boldsymbol{R}\boldsymbol{Q}t}\mathrm{e}^{-\boldsymbol{R}t}
\end{equation}
 so that $\mathrm{e}^{\boldsymbol{R}\boldsymbol{Q}t}=\boldsymbol{\varPhi}(t)\mathrm{e}^{\boldsymbol{R}t}$.
The power series expansion of $\boldsymbol{\varPhi}(t)$ can be written
as $\boldsymbol{\varPhi}(t)=\sum_{n=0}^{\infty}\frac{t^{n}}{n!}\left.\frac{\mathrm{d}^{n}}{\mathrm{d}t^{n}}\boldsymbol{\varPhi}(t)\right|_{t=0}$.
The derivatives are calculated as
\begin{align}
\left.\frac{\mathrm{d}^{n}}{\mathrm{d}t^{n}}\boldsymbol{\varPhi}(t)\right|_{t=0} & =\sum_{k=0}^{n}\binom{n}{k}\left(\boldsymbol{R}-\boldsymbol{R}\boldsymbol{S}\right)^{k}\left(-\boldsymbol{R}\right)^{n-k}\nonumber \\
 & =\sum_{k=0}^{n}\binom{n}{k}\sum_{p_{i}+q_{i}=1}\boldsymbol{R}^{p_{1}}\left(-\boldsymbol{R}\boldsymbol{S}\right)^{q_{1}}\cdots\boldsymbol{R}^{p_{k}}\left(-\boldsymbol{R}\boldsymbol{S}\right)^{q_{k}}\left(-\boldsymbol{R}\right)^{n-k}.\label{eq:phiExp-A}
\end{align}
We can transform the products in (\ref{eq:phiExp-A}) to simpler expressions.
Assume that in the second summation for a fixed $r\in\{0,\ldots,k\}$,
$q_{r}=1$ and $q_{i}=0$ if $i>r$, while the rest of $q_{i}$ are
arbitrary. The sets of $p_{i},q_{i}$ are disjoint for different $r$
values and their union covers all possible $p_{i},q_{i}$ values once.
Using formula (\ref{eq:WApower}) and $\boldsymbol{R}\boldsymbol{S}=\boldsymbol{W}\boldsymbol{A}\boldsymbol{V}$
we find that 
\begin{equation}
\boldsymbol{R}^{p_{1}}\left(-\boldsymbol{R}\boldsymbol{S}\right)^{q_{1}}\cdots\boldsymbol{R}^{p_{k}}\left(-\boldsymbol{R}\boldsymbol{S}\right)^{q_{k}}=(-1)^{\sum q_{i}}\boldsymbol{W}\boldsymbol{A}^{r}\boldsymbol{V}\boldsymbol{R}^{k-r}.
\end{equation}
The sum of all terms corresponding to each $r\neq$0 can be written
as
\begin{equation}
\sum_{l=0}^{r-1}\binom{r-1}{l}(-1)^{l+1}\boldsymbol{W}\boldsymbol{A}^{r}\boldsymbol{V}\boldsymbol{R}^{k-r}=\begin{cases}
-\boldsymbol{W}\boldsymbol{A}\boldsymbol{V}\boldsymbol{R}^{k-1} & \mbox{if}\; r=1\\
0 & \mbox{if}\; r>1
\end{cases}.
\end{equation}
For $r=0$, the sum is $\boldsymbol{R}^{k}$. Therefore the $n$-th
derivative for $n>0$ reads
\begin{align}
\left.\frac{\mathrm{d}^{n}}{\mathrm{d}t^{n}}\boldsymbol{\varPhi}(t)\right|_{t=0} & =\left(-\boldsymbol{R}\right)^{n}+\sum_{k=1}^{n}\binom{n}{k}\left(\boldsymbol{R}^{k}-\boldsymbol{W}\boldsymbol{A}\boldsymbol{V}\boldsymbol{R}^{k-1}\right)\left(-\boldsymbol{R}\right)^{n-k}\\
 & =\sum_{k=0}^{n}\binom{n}{k}\boldsymbol{R}^{k}\left(-\boldsymbol{R}\right)^{n-k}-\sum_{k=1}^{n}\binom{n}{k}\boldsymbol{W}\boldsymbol{A}\boldsymbol{V}\boldsymbol{R}^{k-1}\left(-\boldsymbol{R}\right)^{n-k}\\
 & =(-1)^{n}\boldsymbol{W}\boldsymbol{A}\boldsymbol{V}\boldsymbol{R}^{n-1}.
\end{align}
Substituting the derivatives into the power series we are left with
\begin{equation}
\boldsymbol{\varPhi}(t)=\boldsymbol{I}-\sum_{n=1}^{\infty}\frac{t^{n}}{n!}(-1)^{n-1}\boldsymbol{W}\boldsymbol{A}\boldsymbol{V}\boldsymbol{R}^{n-1}=\boldsymbol{I}-\boldsymbol{W}\boldsymbol{A}\boldsymbol{V}\int_{0}^{t}\mathrm{e}^{-\boldsymbol{R}\tau}\mathrm{d}\tau.
\end{equation}
Multiplying $\boldsymbol{\varPhi}(t)$ from the right by $\mathrm{e}^{\boldsymbol{R}t}$
we get the formula
\begin{equation}
\mathrm{e}^{\boldsymbol{R}\boldsymbol{Q}t}=\mathrm{e}^{\boldsymbol{R}t}-\boldsymbol{R}\boldsymbol{S}\int_{0}^{t}\mathrm{e}^{\boldsymbol{R}(t-\tau)}\mathrm{d}\tau=\mathrm{e}^{\boldsymbol{R}t}-\boldsymbol{R}\boldsymbol{S}\int_{0}^{t}\mathrm{e}^{\boldsymbol{R}\tau}\mathrm{d}\tau.
\end{equation}
With this result the forcing term and the memory kernels become

\begin{gather*}
\boldsymbol{H}(t)\boldsymbol{z}(s)=\boldsymbol{V}\boldsymbol{R}\boldsymbol{Q}\mathrm{e}^{\boldsymbol{R}(t-s)}\boldsymbol{z}(s),\;\boldsymbol{K}(\tau)=\boldsymbol{A},\;\boldsymbol{L}(\tau)=\left(\boldsymbol{V}\mathrm{e}^{\boldsymbol{R}\tau}-\boldsymbol{A}\boldsymbol{V}\int_{0}^{\tau}\mathrm{e}^{\boldsymbol{R}\theta}\mathrm{d}\theta\right)\boldsymbol{v}.
\end{gather*}

\section{Strong continuity of $\mathrm{e}^{\boldsymbol{R}t}$\label{sub:C0semigroup}}

In order to justify our analysis we need to show that $\mathrm{e}^{\boldsymbol{R}t}$
is a strongly continuous semigroup. We use the Hille-Yosida theorem
\citep{Pazy83,HillePhillips57}, which states that $\mathrm{e}^{\boldsymbol{R}t}$
is a strongly continuous semigroup satisfying $\left\Vert \mathrm{e}^{\boldsymbol{R}t}\right\Vert \le M_{0}$
if and only if $\overline{\mathcal{D}(\boldsymbol{R})}=\boldsymbol{Z}$
and 
\begin{equation}
\left\Vert (\boldsymbol{R}-\lambda\boldsymbol{I})^{-n}\right\Vert \le M_{0}\lambda^{-n}\label{eq:resolventCond-HY}
\end{equation}
for $\lambda>0$ and $n=1,2,3,\ldots$, where $\boldsymbol{Z}$ is
defined by equation (\ref{eq:Xspace}).

First we show that if $\boldsymbol{R}$ satisfies (\ref{eq:NegativeRealPart})
and that each eigenvalue of $\boldsymbol{R}$ has finite multiplicity
then condition (\ref{eq:resolventCond-HY}) is satisfied. Using a
linear transformation matrix $\boldsymbol{R}$ can be brought into
its block diagonal Jordan normal form. Each block in the diagonal
of the Jordan normal form corresponds to an eigenvalue $\lambda_{k}$
of $\boldsymbol{R}$ and has size $l_{k}$. The form of such a block
is 
\begin{equation}
\boldsymbol{J}_{k}=\left(\begin{array}{cccc}
\lambda_{k} & 1 & 0 & 0\\
0 & \lambda_{k} & \ddots & 0\\
0 & \ddots & \ddots & 1\\
0 & 0 & 0 & \lambda_{k}
\end{array}\right).
\end{equation}
After inversion the component of $(\boldsymbol{R}-\lambda\boldsymbol{I})^{-n}$
that corresponds to $J_{k}$ becomes
\begin{equation}
\left(\boldsymbol{J}_{k}-\lambda\boldsymbol{I}\right)^{-n}=\left(\begin{array}{cccc}
\left(\lambda_{k}-\lambda\right)^{-n} & -n\left(\lambda_{k}-\lambda\right)^{-n-1} & \cdots & (-1)^{l}\binom{n+l_{k}-1}{l_{k}}\left(\lambda_{k}-\lambda\right)^{-n-l_{k}}\\
0 & \left(\lambda_{k}-\lambda\right)^{-n} & \ddots & \vdots\\
0 & 0 & \ddots & -n\left(\lambda_{k}-\lambda\right)^{-n-1}\\
0 & 0 & 0 & \left(\lambda_{k}-\lambda\right)^{-n}
\end{array}\right).
\end{equation}
The norm of this Jordan block can be estimated by
\begin{equation}
\left\Vert \left(\boldsymbol{J}_{k}-\lambda\boldsymbol{I}\right)^{-n}\right\Vert \le\left|\lambda_{k}-\lambda\right|^{-n}+n\left|\lambda_{k}-\lambda\right|^{-n-1}+\cdots+\binom{n+l_{k}-1}{l_{k}}\left|\lambda_{k}-\lambda\right|^{-n-l_{k}}.
\end{equation}
Note that $\left|\lambda_{k}-\lambda\right|\ge\left|\Re\lambda_{k}-\lambda\right|$.
Since $\lambda>0\ge\Re\lambda_{k}$ one can find an $M_{k}$ such
that 
\begin{equation}
\left\Vert \left(\boldsymbol{J}_{k}-\lambda\boldsymbol{I}\right)^{-n}\right\Vert \le M_{k}\lambda^{-n}
\end{equation}
Considering this estimate for all Jordan blocks we find that 
\begin{equation}
\left\Vert (\boldsymbol{R}-\lambda\boldsymbol{I})^{-n}\right\Vert \le\sup_{k}M_{k}\lambda^{-n}\le M_{0}\lambda^{-n},
\end{equation}
where $M_{0}=\sup_{k}M_{k}$. This proves (\ref{eq:resolventCond-HY}).

To conclude the proof we show that $\overline{\mathcal{D}(\boldsymbol{R})}=\boldsymbol{Z}$.
Again, we use the Jordan normal form. We partition every vector in
$\boldsymbol{Z}$ such that $\boldsymbol{z}=(\boldsymbol{z}_{1},\boldsymbol{z}_{2},\ldots)^{T}$,
where $\boldsymbol{z}_{k}$ are of the size of a Jordan block. Let
$\boldsymbol{z}\in\boldsymbol{Z}$ and construct $\boldsymbol{z}^{P}=(\boldsymbol{z}_{1},\boldsymbol{z}_{2},\ldots,\boldsymbol{z}_{P},0,\ldots)^{T}$
such that it has $P$ number of non-zero components. This guarantees
that for any $\boldsymbol{z}^{P}$, $\left\Vert \boldsymbol{R}\boldsymbol{z}^{P}\right\Vert <\infty$,
hence $\boldsymbol{z}^{P}\in\mathcal{D}(\boldsymbol{R})$. It is also
clear that $\left\Vert \lim_{P\to\infty}\boldsymbol{z}_{P}\right\Vert <\infty$
due to its construction, thus we have shown that $\overline{\mathcal{D}(\boldsymbol{R})}=\boldsymbol{X}$.

\section{Boundedness and smoothness of $\boldsymbol{L}_{j}^{1}(\tau)$\label{sub:L1bound}}

The definition (\ref{eq:L1def}) with (\ref{eq:L0Def}) of $\boldsymbol{L}_{j}^{1}(\tau)$
has two terms both including the expression $\boldsymbol{\Upsilon}_{j}(\tau)=\boldsymbol{V}\int_{0}^{\tau}\mathrm{e}^{\boldsymbol{R}\theta}\boldsymbol{v}_{M+j}\mathrm{d\theta}$.
Therefore we only need to consider $\boldsymbol{\Upsilon}_{j}(\tau)$
in our analysis to show boundedness and smoothness of $\boldsymbol{L}_{j}^{1}(\tau)$.

We use the inverse Laplace transform \citep{Pazy83} to obtain
\begin{equation}
\mathrm{e}^{\boldsymbol{R}\tau}\boldsymbol{x}=\frac{1}{2\pi i}\int_{\Gamma_{0}}\mathrm{e}^{\lambda\tau}\left(\lambda\boldsymbol{I}-\boldsymbol{R}\right)^{-1}\boldsymbol{x}\mathrm{d}\lambda,\label{eq:OpLaplace}
\end{equation}
where $\Gamma_{0}=\{\lambda\in\mathbb{C}:\Re\lambda=\gamma>0\}$.
Using integration rules for the Laplace transform and multiplying
(\ref{eq:OpLaplace}) by vectors $\boldsymbol{v}_{M+j}$ and $\boldsymbol{v}_{\ell}$
from the left and right, respectively, we get
\begin{equation}
\left[\boldsymbol{\Upsilon}_{j}(\tau)\right]_{\ell}=\frac{1}{2\pi i}\int_{\Gamma_{0}}\mathrm{e}^{\lambda\tau}\boldsymbol{v}_{\ell}\cdot\left(\lambda\boldsymbol{I}-\boldsymbol{R}\right)^{-1}\boldsymbol{v}_{M+j}\frac{\mathrm{d}\lambda}{\lambda},\label{eq:projLaplace}
\end{equation}
The equality makes sense if the integral converges. Evaluating the
inverse operator in (\ref{eq:projLaplace}) we find that
\begin{equation}
\left(\lambda\boldsymbol{I}-\boldsymbol{R}\right)^{-1}\left(\begin{array}{c}
\boldsymbol{0}\\
\boldsymbol{n}_{j}
\end{array}\right)=\left(\begin{array}{c}
-(\boldsymbol{K}+\lambda\boldsymbol{C}+\lambda^{2})^{-1}\boldsymbol{n}_{j}\\
-\lambda(\boldsymbol{K}+\lambda\boldsymbol{C}+\lambda^{2})^{-1}\boldsymbol{n}_{j}
\end{array}\right).\label{eq:ResolvForm}
\end{equation}
Note that it is sufficient to consider $\left[\boldsymbol{\Upsilon}_{j}(\tau)\right]_{M+l}$
since $\left[\boldsymbol{\Upsilon}_{j}(\tau)\right]_{l}$ is the integral
of $\left[\boldsymbol{\Upsilon}_{j}(\tau)\right]_{M+l}$. Substituting
(\ref{eq:ResolvForm}) into (\ref{eq:projLaplace}) we are left with
\begin{align}
\left[\boldsymbol{\Upsilon}_{j}(\tau)\right]_{M+l} & =\frac{1}{2\pi i}\int_{\Gamma_{0}}\mathrm{e}^{\lambda\tau}\boldsymbol{n}_{l}\cdot(\boldsymbol{K}+\lambda\boldsymbol{C}+\lambda^{2})^{-1}\boldsymbol{n}_{j}\mathrm{d}\lambda.
\end{align}
The integral of the inverse Laplace transform (\ref{eq:projLaplace})
converges for all $t\ge0$ if 
\begin{equation}
\left|\lambda\boldsymbol{n}_{l}\cdot(\boldsymbol{K}+\lambda\boldsymbol{C}+\lambda^{2})^{-1}\boldsymbol{n}_{j}\right|\le M_{j,l}\;\mbox{for}\;\lambda\in\Gamma_{0}
\end{equation}
and for $j,l=1,\ldots,M$. Indeed, by estimating the bound we get
\begin{equation}
\left|\int_{0}^{\tau}\left[\boldsymbol{\Upsilon}_{j}(\tau)\right]_{M+l}\right|\le\frac{1}{2\pi i}M_{j,l}\left|\int_{\Gamma_{0}}\mathrm{e}^{\lambda\tau}\frac{\mathrm{d}\lambda}{\lambda}\right|=M_{j,l}.
\end{equation}
This implies that $\boldsymbol{L}_{j}^{1}$ is bounded. 

If the stronger conditions (\ref{eq:SectorialCond}) and (\ref{eq:SmoothCond})
are satisfied, we can alter the contour of integration so that it
is within the left half of the complex plane, $\Gamma_{\delta}=\{\lambda\in\mathbb{C}:\left|\arg\lambda\right|=\nicefrac{\pi}{2}+\delta\}$,
where $\delta>0$ is sufficiently small so that all the eigenvalues
of $\boldsymbol{R}$ are on the left of $\Gamma_{\delta}$. We parametrise
the contour $\Gamma_{\delta}$ by $\lambda=\kappa\mathrm{e}^{\pm i(\nicefrac{\pi}{2}+\delta)}$,
$\kappa\ge0$. Using this contour we find that the derivative
\begin{equation}
\left[\boldsymbol{\Upsilon}_{j}(\tau)\right]_{M+l}=\frac{1}{2\pi i}\int_{\Gamma_{\delta}}\mathrm{e}^{\lambda\tau}\lambda\boldsymbol{n}_{l}\cdot(\boldsymbol{K}+\lambda\boldsymbol{C}+\lambda^{2})^{-1}\boldsymbol{n}_{j}\mathrm{d}\lambda
\end{equation}
can be estimated by
\begin{equation}
\left|\left[\boldsymbol{\Upsilon}_{j}(\tau)\right]_{M+l}\right|\le\frac{1}{\pi}\int_{0}^{\infty}M_{j,l}\mathrm{e}^{-\kappa\tau\cos2\delta}\mathrm{d}\kappa=\frac{1}{\tau}\left(\frac{M_{j,l}}{\pi\cos2\delta}\right).
\end{equation}
This means that the derivative of $\boldsymbol{L}_{j}^{1}$ is also
bounded but only for $\tau>0$. Due to $\mathrm{e}^{\boldsymbol{R}\tau}$
being a strongly continuous semigroup (Theorem 2.5.2 in \citet{Pazy83})
higher derivatives are also bounded 
\begin{equation}
\frac{1}{n!}\frac{\mathrm{d}^{n-1}}{\mathrm{d}\tau^{n-1}}\left|\left[\boldsymbol{\Upsilon}_{j}(\tau)\right]_{M+l}\right|\le\left(\frac{C_{j,l}\mathrm{e}}{\tau}\right)^{n},
\end{equation}
for $C_{j,l}<\infty$ which proves that $\boldsymbol{L}_{j}^{1}$
is smooth for $\tau>0$.

\section{Discontinuity of the memory kernel of the pre-tensed string\label{sec:L1Lim}}

Here we show that the memory kernel $\boldsymbol{L}^{1}(\tau)$ of
the bowed string is discontinuous at $\tau=0$. In equation (\ref{eq:StringL2intGen})
the terms that cause discontinuity are divided by the lowest power
of $\omega_{k}$. The other terms are continuous and add up to zero
at $\tau=0$. Therefore, after using $D_{k}=D$ and $\omega_{k}=ck\pi$
the following identity holds: 
\begin{equation}
\lim_{\tau\to0+}\left[\boldsymbol{L}^{1}(\tau)\right]_{2}=\lim_{\tau\to0+}\sum_{k=1}^{\infty}\kappa_{k}(\tau),\;\kappa_{k}(\tau)=\frac{e^{-k\pi cD\tau}\sin^{2}k\pi\xi^{\star}}{k\pi c\sqrt{1-D^{2}}}\sin\left(k\pi c\sqrt{1-D^{2}}\tau\right).\label{eq:C-L1lim}
\end{equation}
One can expand $\kappa_{k}(\tau)$ in (\ref{eq:C-L1lim}) as a sum
of exponentials
\begin{multline}
\kappa_{k}(\tau)=\frac{1}{4k\pi c\sqrt{D^{2}-1}}\biggl(\mathrm{e}^{k\pi c\left(\sqrt{D^{2}-1}-D\right)\tau}-\mathrm{e}^{-k\pi c\left(\sqrt{D^{2}-1}+D\right)\tau}+\frac{1}{2}\mathrm{e}^{-k\pi\left(c\left(\sqrt{D^{2}-1}+D\right)\tau+2i\zeta\right)}\\
-\frac{1}{2}\mathrm{e}^{k\pi\left(c\left(\sqrt{D^{2}-1}-D\right)\tau-2i\zeta\right)}+\frac{1}{2}\mathrm{e}^{-k\pi\left(c\left(\sqrt{D^{2}-1}+D\right)\tau-2i\zeta\right)}-\frac{1}{2}\mathrm{e}^{k\pi\left(c\left(\sqrt{D^{2}-1}-D\right)\tau+2i\zeta\right)}\biggr).
\end{multline}
Since $\sum_{k=1}^{\infty}\frac{\mathrm{e}^{ka}}{k}=-\log(1-\mathrm{e}^{a})$,
the limit (\ref{eq:C-L1lim}) can be written as
\begin{equation}
\lim_{\tau\to0+}\left[\boldsymbol{L}^{1}(\tau)\right]_{2}=\lim_{\tau\to0+}\frac{-1}{4\pi c\sqrt{D^{2}-1}}\left\{ \log(1-\mathrm{e}^{a_{1}})-\log(1-\mathrm{e}^{a_{2}})+\frac{1}{2}\sum_{l=3}^{6}(-1)^{l+1}\log(1-\mathrm{e}^{a_{l}})\right\} ,
\end{equation}
where
\begin{eqnarray}
a_{1}=\pi c\left(\sqrt{D^{2}-1}-D\right)\tau, & a_{2}=-\pi c\left(\sqrt{D^{2}-1}+D\right)\tau,\\
a_{3,4}=\mp\pi\left(c\left(\sqrt{D^{2}-1}+D\right)\tau+2i\zeta\right), & a_{5,6}=\mp\pi\left(c\left(\sqrt{D^{2}-1}+D\right)\tau-2i\zeta\right),
\end{eqnarray}
Discontinuity occurs if the path of $1-\mathrm{e}^{a_{l}(\tau)}$
crosses the non-positive real axis (including zero) at $\tau=0$.
This is possible for $a_{1}$ and $a_{2}$ only if $0<\xi<1$, $\xi\neq\nicefrac{1}{2}$.
Since we are taking a limit, it is sufficient to use a first order
approximation at $\tau=0$, that is, $1-\mathrm{e}^{a_{l}(\tau)}\approx-a_{l}(\tau)$.
Also note that $\log x=\log\left|x\right|+i\arg x$ and that $\left|a_{1}\right|=\left|a_{2}\right|=\pi c\tau$.
The limit therefore becomes
\begin{equation}
\lim_{\tau\to0+}\left[\boldsymbol{L}^{1}(\tau)\right]_{2}=\frac{-1}{4\pi c\sqrt{1-D^{2}}}\left(\arg\left(\sqrt{D^{2}-1}-D\right)-\arg\left(-\sqrt{D^{2}-1}-D\right)\right).
\end{equation}
Assuming that $D=\cos\phi$, $0\le\phi\le\nicefrac{\pi}{2}$, we get
$\arg\left(\sqrt{D^{2}-1}-D\right)=\pi-\phi$ and $\arg\left(-\sqrt{D^{2}-1}-D\right)=\pi+\phi$,
hence
\begin{equation}
\lim_{\tau\to0+}\left[\boldsymbol{L}^{1}(t)\right]_{2}=\frac{\cos^{-1}D}{2\pi c\sqrt{1-D^{2}}}.
\end{equation}

\bibliographystyle{rspublicnat}
\bibliography{AllRef}

\end{document}